\documentclass{amsart}
\usepackage{amsmath,amssymb,amsthm,amsxtra,amscd}

\DeclareMathOperator{\A}{\mathbf{A}}

\DeclareMathOperator{\B}{\mathbf{B}}

         \newcommand{\Btimes}{\B_{\times}}

\DeclareMathOperator{\bsv}{\boldsymbol{v}}

\DeclareMathOperator{\C}{\mathbf{C}}

\DeclareMathOperator{\co}{\mathrm{co}}
\DeclareMathOperator{\coim}{coim}
\DeclareMathOperator{\coker}{coker}

\DeclareMathOperator{\D}{\mathbf{D}}

\DeclareMathOperator{\dist}{dist}
\DeclareMathOperator{\Dw}{\mathbf{D}^{\boldsymbol{w}}}
\DeclareMathOperator{\E}{\mathbf{E}}
\DeclareMathOperator{\ebold}{\mathbf{e}}
         
\DeclareMathOperator{\eB}{\mathbf{eB}}

\DeclareMathOperator{\eUb}{\mathbf{eU}^{\mathrm{b}}}
\DeclareMathOperator{\EGamma}{\boldsymbol{E}\boldsymbol{\Gamma}}

\DeclareMathOperator{\End}{End}

 \DeclareMathOperator{\fil}{fil}
\DeclareMathOperator{\F}{\mathbf{F}}

\DeclareMathOperator{\Fun}{Fun}
\DeclareMathOperator{\G}{\mathit{G}}

         \newcommand{\Gnc}{G^{-\infty}}

         \newcommand{\GncGammazero}{\G^{-\infty}_{\Gamma,0}}

         \newcommand{\Gncm}{\G^{-\infty}_{m}}
         \newcommand{\Gnctimes}{\G^{-\infty}_{\times}}

\DeclareMathOperator{\h}{\mathit{h}}
\DeclareMathOperator{\hlf}{\h^{\textit{lf}}}

\DeclareMathOperator{\Hom}{Hom}

\DeclareMathOperator{\id}{id}

\DeclareMathOperator{\im}{im}

\DeclareMathOperator{\incl}{incl}

\DeclareMathOperator{\K}{\mathit{K}}
         \newcommand{\Knc}{\K^{-\infty}}

         \newcommand{\KncGammazero}{\K^{-\infty}_{\Gamma,0}}

\DeclareMathOperator{\mB}{\mathbf{mB}}
\DeclareMathOperator{\mUb}{\mathbf{mU}^{\mathrm{b}}}
         \newcommand{\mbold}{\mathbf{m}}

\DeclareMathOperator{\Mod}{\mathbf{Mod}}
\DeclareMathOperator{\Modf}{\mathbf{Modf}}

\DeclareMathOperator{\point}{point}

\DeclareMathOperator{\pr}{pr}
         
\DeclareMathOperator{\proj}{proj}

\DeclareMathOperator{\Spt}{Spt} 

         \newcommand{\subdot}{\boldsymbol{\cdot}}
\DeclareMathOperator{\supp}{supp}

\DeclareMathOperator{\U}{\mathbf{U}}

\DeclareMathOperator{\Ub}{\mathbf{U}^{\mathrm{b}}}

\DeclareMathOperator{\vB}{\mathbf{\boldsymbol{v}B}}

\DeclareMathOperator{\w}{\boldsymbol{w}}

\DeclareMathOperator{\vD}{\mathbf{\boldsymbol{v}D}}
\DeclareMathOperator{\vDw}{\mathbf{\boldsymbol{v}D^{\boldsymbol{w}}}}
\DeclareMathOperator{\vE}{\mathbf{\boldsymbol{v}E}}

         \newcommand{\bfw}{\boldsymbol{w}}

\DeclareMathOperator{\wB}{\mathbf{\boldsymbol{w}B}}
\DeclareMathOperator{\wD}{\mathbf{\boldsymbol{w}D}}

\DeclareMathOperator{\Z}{\mathbf{Z}}

\DeclareMathOperator*{\one}{1}
\newcommand{\onehatplace}[1]%           % "frown" over 1
{ \one^{\substack{#1 \\ \frown}} }

\DeclareMathOperator*{\bones}{\times}
\newcommand{\undertimes}[1]%            % X with stuff under it
{ \bones_{#1} }

\DeclareMathOperator*{\bowl}{\cup}
\newcommand{\undercup}[1]%              % small U with stuff under it
{ \bowl_{#1} }

\DeclareMathOperator*{\arch}{\cap}
\newcommand{\undercap}[1]%              % small "cap" with stuff under it
{ \arch_{#1} }

\newcommand{\pull}%                              % extension for arrows
{\!\!\! -\!\!\! -\!\!\! -\!\!\!}

\DeclareMathOperator*{\holimprep}{holim}                       % holim
\newcommand{\holim}[1]%
{\displaystyle\holimprep_{\substack{\leftarrow \pull - \\ #1}} \,
}

\DeclareMathOperator*{\hocolimprep}{hocolim}                   % hocolim
\newcommand{\hocolim}[1]%
{\displaystyle\hocolimprep_{\substack{- \pull \rightarrow \\ #1}}
\, }

\DeclareMathOperator*{\plainlim}{lim}                           % lim
\newcommand{\contralim}[1]%
{\displaystyle\plainlim_{\substack{\leftarrow \pull - \\ #1}} \, }

\DeclareMathOperator*{\plaincolim}{colim}                       % colim
\newcommand{\colim}[1]%
{\displaystyle\plaincolim_{\substack{- \pull \rightarrow \\ #1}}
\, }

\DeclareMathOperator*{\laxlimplain}{laxlim}                     % laxlim
\newcommand{\laxlim}[1]%
{\displaystyle\laxlimplain_{\substack{\leftarrow \pull - \\ #1}}
\, }

\providecommand{\bysame}{\makebox[3em]{\hrulefill}\thinspace}

\theoremstyle{plain}
\newtheorem{Thm}{Theorem}[section]
\newtheorem*{ThmPlain}{Theorem}

\newtheorem{Cor}[Thm]{Corollary}

\newtheorem{Lem}[Thm]{Lemma}
\newtheorem{Prop}[Thm]{Proposition}

\theoremstyle{definition}
\newtheorem{Def}[Thm]{Definition}

\newtheorem{Ex}[Thm]{Example}

\newtheorem{Rem}[Thm]{Remark}

\newtheorem{SatAx}[Thm]{Saturation axiom}
\newtheorem{ExtAx}[Thm]{Extension axiom}
\newtheorem{CylAx}[Thm]{Cylinder axiom}

\theoremstyle{remark}
\newtheorem{Not}[Thm]{Notation}

\newtheoremstyle{freestylethm}{6pt}{6pt}{\itshape}{}%
                {\bfseries}{}{.5em}{\thmnote{#3}}
\theoremstyle{freestylethm}

\newcommand{\define}[1]{\textit{#1}}

\newif\ifShowLabels
\ShowLabelstrue
\newcommand{\TeXref}[1]{
\marginpar{\scriptsize \texttt{#1}}}

\newcommand{\SecRef}[2]{\section{#1}\label{S:#2}%
\ifShowLabels \TeXref{{S:#2}} \fi}

\newcommand{\refT}[1]{\textup{\ref{T:#1}}}
\newcommand{\refL}[1]{\textup{\ref{L:#1}}}
\newcommand{\refD}[1]{\textup{\ref{D:#1}}}
\newcommand{\refC}[1]{\textup{\ref{C:#1}}}
\newcommand{\refE}[1]{\textup{\ref{E:#1}}}
\newcommand{\refP}[1]{\textup{\ref{P:#1}}}
\newcommand{\refR}[1]{\textup{\ref{R:#1}}}

\newenvironment{ThmRef}[1]%
{ \begin{Thm} \label{T:#1}
\ifShowLabels \TeXref{T:#1} \fi }%
{ \end{Thm} }
\newenvironment{DefRef}[1]%
{ \begin{Def} \label{D:#1}
\ifShowLabels \TeXref{D:#1} \fi }%
{ \end{Def} }
\newenvironment{LemRef}[1]%
{ \begin{Lem} \label{L:#1}
\ifShowLabels \TeXref{L:#1} \fi }%
{ \end{Lem} }
\newenvironment{CorRef}[1]%
{ \begin{Cor} \label{C:#1}
\ifShowLabels \TeXref{C:#1} \fi }%
{ \end{Cor} }
\newenvironment{RemRef}[1]%
{ \begin{Rem} \label{R:#1}
\ifShowLabels \TeXref{R:#1} \fi }%
{ \end{Rem} }
\newenvironment{PropRef}[1]%
{ \begin{Prop} \label{P:#1}
\ifShowLabels \TeXref{P:#1} \fi }%
{ \end{Prop} }
\newenvironment{ExRef}[1]%
{ \begin{Ex} \label{E:#1}
\ifShowLabels \TeXref{E:#1} \fi  }%
{ \end{Ex} }
\newenvironment{NotRef}[1]%
{ \begin{Not} \label{N:#1}
\ifShowLabels \TeXref{N:#1} \fi }%
{ \end{Not} }

\newenvironment{ThmRefName}[2]%
{ \begin{Thm} [#2]\label{T:#1}
\ifShowLabels \TeXref{T:#1} \fi }%
{ \end{Thm} }
\newenvironment{DefRefName}[2]%
{ \begin{Def} [#2]\label{D:#1}
\ifShowLabels \TeXref{D:#1} \fi }%
{ \end{Def} }
\newenvironment{LemRefName}[2]%
{ \begin{Lem} [#2]\label{L:#1}
\ifShowLabels \TeXref{L:#1} \fi }%
{ \end{Lem} }
\newenvironment{CorRefName}[2]%
{ \begin{Cor} [#2]\label{C:#1}
\ifShowLabels \TeXref{C:#1} \fi }%
{ \end{Cor} }
{ \begin{Rem} [#2]\label{R:#1}
\ifShowLabels \TeXref{R:#1} \fi }%
{ \end{Rem} }
{ \begin{Prop} [#2]\label{P:#1}
\ifShowLabels \TeXref{P:#1} \fi }%
{ \end{Prop} }
{ \begin{Ex} [#2]\label{E:#1}
\ifShowLabels \TeXref{E:#1} \fi }%
{ \end{Ex} }

\ShowLabelsfalse

\begin{document}

\title[Controlled algebraic \textit{G}-theory, I]
{Controlled algebraic \textit{G}-theory, I}
\author{}
\author[Gunnar Carlsson]{Gunnar Carlsson}
\address{Department of Mathematics\\
Stanford University\\ Stanford\\ CA 94305}
\email{gunnar@math.stanford.edu}
\author[Boris Goldfarb]{Boris Goldfarb}
\address{Department of Mathematics and Statistics\\
SUNY\\ Albany\\ NY 12222} \email{goldfarb@math.albany.edu}
\date{\today}
\thanks{The authors acknowledge support from the
National Science Foundation.}

\begin{abstract}
This paper extends the notion of geometric control in algebraic
$K$-theory from additive categories with split exact sequences to
other exact structures. In particular, we construct exact categories
of modules over a Noetherian ring filtered by subsets of a metric
space and sensitive to the large scale properties of the space.  The
algebraic $K$-theory of these categories is related to the
bounded $K$-theory of geometric modules of Pedersen and Weibel the way $G$-theory is
classically related to $K$-theory.  We recover familiar results in
the new setting, including the nonconnective bounded excision and
equivariant properties.
We apply the results to the $G$-theoretic Novikov conjecture which is
shown to be stronger than the usual $K$-theoretic conjecture.
\end{abstract}

\maketitle

\tableofcontents

\SecRef{Introduction}{Intro}

Since the invention of algebraic $K$-groups of a ring defined
using the finitely generated projective $R$-modules, there existed
a companion $K$-theory defined using arbitrary finitely generated
$R$-modules, called $G$-theory. Its usefulness comes from the
computational tool available in $G$-theory, the localization exact
sequence, and the close relation to $K$-theory via the Cartan map
which becomes an isomorphism when $R$ is a regular ring. The
recent success of controlled $K$-theory in algebra and topology,
where the ring involved is usually the regular ring of integers
$\mathbb{Z}$, makes it natural to look for a similar controlled
analogue of $G$-theory. This paper constructs and exploits such an
analogue.

The bounded control is introduced by fixing a basis $B$ in a free
module $M$ and defining a locally finite set function $s \colon B
\to X$ into a metric space $X$. The control comes from
restrictions on the maps one allows between the based modules.
Since each element $x$ in $M$ is written uniquely as a sum $x =
\sum_{b \in B} r_b b$, there is the notion of support, $\supp (x)$,
which is the set of all points
$s (b)$ in $X$ with $b \in B$ such that $r_b \ne 0$.

For two sets of choices $(M_i, B_i, s_i)$, $i=1$, $2$, an
$R$-homomorphism $\phi \colon M_1 \to M_2$ is \define{bounded} if
there is a number $D > 0$ such that for every $b \in B_1$ the
support $\supp s_2 (\phi (b))$ is contained in the metric ball of
radius $D$ centered at $s_1 (b)$.

The triples $(M, B, s)$
and the bounded homomorphisms form the \define{bounded
category} $\mathcal{B} (X,R)$. It is in fact an additive category
since the direct products can be defined in the evident way. To
each additive category $\mathcal{A}$, one associates a
sequence of groups $K_i (\mathcal{A})$, $i \in \mathbb{Z}$, or
rather a nonconnective spectrum $\Spt (\mathcal{A})$ whose stable
homotopy groups are $K_i (\mathcal{A})$, as in \cite{ePcW:89}.
This construction applied to the bounded category $\mathcal{B} (X,R)$ gives the
\define{bounded algebraic $K$-theory} $K_i (X,R)$.

The general goal of this paper is to construct larger categories
associated to a metric space $X$ and a Noetherian ring $R$
and to recover in this context the basic results from bounded
$K$-theory. We are mostly concerned with controlled
excision established in section 3.
In many ways these categories are more flexible than the bounded
categories and allow application of recent powerful results in
algebraic $K$-theory.
Their properties are essential for our study of the Borel isomorphism conjecture
continued elsewhere but indicated in section 4.
In the same section, we prove the integral Novikov
conjecture in this context.
It turns out that earlier results asserting the split injectivity of the assembly map follow
from this result, due to the fact that the natural transformation from $K$-theory to $G$-theory is an equivalence for a regular Noetherian ring.

The following is a sketch of the constructions and results of the paper.

First notice that, given a triple $(M, B, s)$ in $\mathcal{B}
(X,R)$, to every subset $S \subset X$ there is associated a free
submodule $M (S)$ generated by those $b \in B$ with the property
$s(b) \in S$. The restriction to bounded homomorphisms can be
described entirely in terms of these submodules. We generalize
this as follows. The objects of the new category $\Ub (X,R)$ are
left $R$-modules $M$ filtered by the subsets of $X$ in the sense
that they are functors from the category of subsets of $X$ and inclusions
to the category of submodules of $M$ and inclusions
for which the value on the whole space $X$ is the whole module $M$ and the value on the empty set $\varnothing$ is the zero submodule.  By abuse of notation
we usually denote the functor by the same letter $M$.  We also
make several additional assumptions spelled out in Definition
\refD{RealBCprelim}, in particular, that the values on the bounded subsets are finitely generated submodules.

The morphisms in $\Ub (X,R)$ are the
$R$-homomorphisms $\phi \colon M_1 \to M_2$ for which there exists a number $D \ge 0$
such that the image
$\phi (M_1 (S))$ is contained in the submodule $M_2 (S [D])$ for
all subsets $S \subset X$. Here $S[D]$ stands for the
metric $D$-enlargement of $S$ in $X$. In this context we say a submodule $N \subset M$
is \textit{supported on a subset} $S \subset X$ if $N \subset M(S)$.

The
\define{boundedly controlled} category $\B (X, R)$ is the full
subcategory of $\Ub (X, R)$ on filtered modules $M$ generated by
elements supported on subsets of diameter less than $d$ for some
number $d > 0$ specific to $M$.

The additive structure on $\B (X, R)$ gives it the \textit{split
exact structure} where the admissible monomorphisms are all split
monics and admissible epimorphisms are all split epis.  In order
to describe a different Quillen exact structure on $\B (X, R)$, we
define an additional property a boundedly controlled homomorphism
$\phi \colon M_1 \to M_2$ in $\Ub (X, R)$ may or may not have:
$\phi$ is
\define{boundedly bicontrolled} if there exists a number $D \ge 0$
such that
\[
\phi (M_1 (S)) \subset M_2 (S [D])
\]
and
\[
\phi (M_1) \cap M_2 (S) \subset \phi M_1 (S[D])
\]
for all subsets $S$ of $X$.
We define the admissible monomorphisms in the new Quillen exact structure to be the
boundedly bicontrolled
injections of modules, both for the case of
$\B (X,R)$ and $\Ub (X,R)$.
We define the admissible epimorphisms in $\Ub (X,R)$
to be the boundedly bicontrolled surjections.
The admissible
epimorphisms in $\B (X,R)$ are the boundedly bicontrolled
surjections with kernels in $\B (X,R)$.
In both cases the exact
sequences are simply the short exact sequences when viewed as
sequences in $\Ub (X,R)$ so that all kernels and cokernels are
well-defined filtered submodules in the respective category.
Notice that split injections and surjections are boundedly
bicontrolled, so the split exact structure is an exact subcategory
of the new one.

Recall that a map $f \colon X \to Y$
of metric spaces is
\define{bi-Lipschitz} if there is a number $k \ge 1$ such that
\[
k^{-1} \dist(x_1, x_2) \le \dist (f(x_1), f(x_2)) \le k \dist(x_1,
x_2)
\]
for all $x_1$, $x_2 \in X$.
More generally, $f$ is \define{quasi-bi-Lipschitz} if there is a real positive function $l$ such that
\begin{gather}
\dist(x_1, x_2) \le r \ \Longrightarrow
\dist (f(x_1), f(x_2)) \le l(r), \notag \\
\dist (f(x_1), f(x_2)) \le r \ \Longrightarrow \ \dist(x_1, x_2)
\le l(r). \notag
\end{gather}
For example, a bounded function $f \colon X \to X$, with the property
$\dist(x, f(x)) \le D$ for all $x \in X$ and a fixed number $D \ge 0$, is
quasi-bi-Lipschitz with $l(r)=r+2D$. An isometry $g \colon X \to Y$
is quasi-bi-Lipschitz with $l (r) = r$.

Both constructions, $\Ub (X,R)$ and $\B (X,R)$, are functorial in the metric space variable
$X$ with respect to quasi-bi-Lipschitz
maps, as should be expected from \cite{ePcW:89}.
Recall that an exact functor between Quillen exact categories is an additive functor which sends exact sequences to exact sequences.
So to each quasi-bi-Lipschitz map $f \colon X \to Y$ one associates an exact functor $f_{\ast} \colon \B (X,R) \to \B (Y,R)$.
For example, if $Z$ is a metric subspace of $X$ then the isometric inclusion $i \colon Z \to X$ induces
an exact functor $i_{\ast} \colon \B (Z,R) \to \B (X,R)$.

To an exact category $\E$, one associates a sequence of groups
$K_i (\E)$, $i \ge 0$, as in Quillen \cite{dQ:73} or a connective
spectrum $K (\E)$ whose stable homotopy groups are $K_i (\E)$.  If
the exact structure is split, these groups are the same as the
$K$-groups of $\E$ as an additive category.
When this construction is applied to the exact category $\B (X,R)$, we call the result the \textit{connective controlled G-theory} of $X$ and denote the spectrum by $G (X, R)$.

A proper metric space is a metric space where all closed
metric balls in are compact.
Let $X$ be a proper metric space.
Suppose $Z$ is a metric subspace of $X$.  There is a construction of
an exact category ${\B}/{\Z}$ associated to $Z$ and an exact functor $\B
(X,R) \to {\B}/{\Z}$ such that the following is true.

\begin{ThmPlain}[Localization, Corollary \refC{LocArb}] The sequence
\[
G (Z, R) \longrightarrow G (X, R) \longrightarrow K ({\B}/{\Z})
\]
is a homotopy fibration.
\end{ThmPlain}

The Localization Theorem can be used to construct nonconnective deloopings of $G
(X,R)$.
We will indicate the
corresponding nonconnective spectra with superscripts ``$-\infty$''.
The construction is similar to the $K$-theory delooping using bounded $K$-theory due to Pedersen--Weibel \cite{ePcW:85}.
Therefore, there is a natural transformation $\Knc (X,R) \to \Gnc
(X,R)$.

The following is the analogue of a major tool in many proofs of the
Novikov conjecture.
If a proper metric space $X$ is the union of
proper metric subspaces $X_1$ and $X_2$, let $\B (X_1,X_2;R)$ stand for the full
subcategory of $\B (X,R)$ on the modules supported on the
intersection of bounded enlargements of $X_1$ and $X_2$ and let $G
(X_1,X_2;R)$ denote its $K$-theory.

\begin{ThmPlain}[Nonconnective controlled excision, Theorem \refT{Exc}]
There is a homotopy pushout
\[
\begin{CD}
\Gnc (X_1, X_2; R) @>>> \Gnc (X_1, R) \\
@VVV @VVV \\
\Gnc (X_2, R) @>>> \Gnc (X, R)
\end{CD}
\]
\end{ThmPlain}

Finally, we describe the application to splitting integral $G$-theoretic assembly maps.
There is a close relation to the same problem in $K$-theory.

Earlier applications of bounded $K$-theory to conjectures of Novikov type
use in a critical way the existence of equivariant versions of the bounded
$K$-theory functors attached to actions of discrete groups of isometries.
The applications we envision of the present theory will also require such
a theory, and we develop it in the last section of the paper. It turns out
that we will need to develop the equivariant theory for a more general
class of actions than isometric actions, namely the class of actions by
discrete groups on metric spaces by quasi-bi-Lipschitz equivalences.  In
carrying this out, we find that we obtain a novel exact structure $\B (R[\Gamma])$ on the a
category of (not necessarily projective) finitely generated modules over a
group ring $R[\Gamma]$, where $R$ is a Noetherian ring and $\Gamma$ is a
discrete group.

Recall that
the integral assembly map in algebraic $K$-theory
\[
A_K \colon B\Gamma_{+} \wedge \Knc (R) \longrightarrow \Knc
(R[\Gamma])
\]
is defined for any group $\Gamma$ and any ring $R$ and relates the
homology of $\Gamma$ with coefficients in the $K$-theory of $R$ to
the $K$-theory of the group ring. The \define{integral Novikov
conjecture} for $\Gamma$ is the statement that $A_K$ is a split
injection of spectra. It is speculated to be true whenever $\Gamma$
is a discrete torsion-free group.

For a Noetherian ring $R$
and the spectrum $\Gnc (R[\Gamma])$ defined as $\Knc \B (R[\Gamma])$
there is a similar map
\[
A_G \colon B\Gamma_{+} \wedge \Gnc (R) \longrightarrow \Gnc
(R[\Gamma])
\]
which we call the \define{integral assembly map} in algebraic
$G$-theory.  In this paper we show that it is a split injection
for many geometric groups.

\begin{ThmPlain}
Let $\Gamma$ be a discrete group of finite asymptotic dimension
and a finite classifying space.  Let $R$ be a Noetherian ring.
Then the assembly map $A_G$ is a split injection.
\end{ThmPlain}

It turns out that earlier results asserting the split injectivity of the
assembly map follow from this result, due to the fact that the natural
transformation from $K$-theory to $G$-theory for rings is an
equivalence when the ring is regular Noetherian, for example the integers $\mathbb{Z}$.

Indeed, notice that from the commutative
square
\[
\begin{CD}
B\Gamma_{+} \wedge \Knc (R) @>{A_K}>>
\Knc (R[\Gamma]) \\
@V{\simeq}VV @VV{C}V \\
B\Gamma_{+} \wedge \Gnc (R) @>{A_G}>> \Gnc (R[\Gamma])
\end{CD}
\]
the assembly map $A_G$ is, up to homotopy, the composition of
$A_K$ followed by the Cartan map
\[
C \colon \Knc (R[\Gamma]) \longrightarrow \Gnc
(R[\Gamma]).
\]
If $A_G$ is a split injection, it follows that $A_K$
is a split injection.

We are grateful to Marco
Schlichting for showing us the preliminary version of his thesis
\cite{mS:99} and several very fruitful discussions.
We are grateful for the critique and suggestions of the referee and the editor which have improved the paper.

The authors acknowledge support from the
National Science Foundation.

\SecRef{Controlled categories of filtered objects}{OthObj}

This work is motivated by the delooping of algebraic $K$-theory of
a small additive category in \cite{ePcW:85} and, in particular,
the introduction of bounded control in a cocomplete additive
category $\A$ which we briefly recall.

A category is \textit{cocomplete} if it contains colimits of
arbitrary small diagrams, cf.~Mac Lane \cite{sM:71}, chapter V.

\begin{DefRefName}{EPCW}{Pedersen--Weibel}
Let $X$ be a proper metric space, in the sense that all closed
metric balls in $X$ are compact.
An $X$-graded object is a function $F$ from the set $X$ to the set of
objects of $\A$ such that the set $\{x \in S \mid F(x) \ne 0 \}$ is finite for
every bounded $S \subset X$.
We will also refer to the object
\[
F = \bigoplus_{x \in X} F(x)
\]
in $\A$ as an $X$-graded object.

The $X$-graded objects form a new category $\mathcal{B} (X,\A)$.
The morphisms
are collections of $\A$-morphisms $f (x,y) \colon F(x) \to G(y)$
with the property that there is a number $D>0$ such that $f(x,y) =
0$ if $\dist (x,y) >D$.

If $\B$ is a subcategory of $\A$ closed
under the direct sum, one obtains the additive bounded category
$\mathcal{B} (X,\B)$ as the full subcategory of $\mathcal{B} (X,\A)$
on objects $F$ with $F(x) \in \B$ for all $x \in X$.  Notice that
$\B$ does not need to be cocomplete.

The \textit{bounded algebraic
$K$-theory} $K (X,\B)$ is the $K$-theory spectrum associated to
the additive category
$\mathcal{B} (X,\B)$.
\end{DefRefName}

To generalize this construction from additive to more general exact
categories $\E$, first notice the following.  Given an object $F$ in
$\mathcal{B} (X,\B)$, to every subset $S \subset X$ there is
associated a direct sum
\[
F(S) = \bigoplus_{x
\in S} F(x).
\]
Since the condition $f(x,y) =
0$ if $\dist (x,y) >D$ is equivalent to the condition that $f (F(x)) \subset F(x[D])$ or,
more generally,
\[
f(F(S)) \subset F(S[D]),
\]
the restriction from arbitrary to bounded morphisms can be described
entirely in terms of the subobjects $F(S)$.

We start by recalling definitions and several standard facts about exact and abelian categories.

\begin{DefRefName}{QEC}{Quillen exact categories}
Let $\C$ be an additive category. Suppose $\C$ has two classes of
morphisms $\mbold (\C)$, called \define{admissible monomorphisms},
and $\ebold (\C)$, called \define{admissible epimorphisms}, and a
class $\mathcal{E}$ of \define{exact} sequences, or extensions, of
the form
\[
C^{\subdot} \colon \quad C' \xrightarrow{i} C \xrightarrow{j} C''
\]
with $i \in \mbold (\C)$ and $j \in \ebold (\C)$ which satisfy the
three axioms:
\begin{enumerate}
\item[a)] any sequence in $\C$ isomorphic to a sequence in $\mathcal{E}$
is in $\mathcal{E}$; the canonical sequence
\[
\begin{CD}
C' @>{\incl_1}>> C' \oplus C'' @>{\proj_2}>> C''
\end{CD}
\]
is in $\mathcal{E}$; for any sequence $C^{\subdot}$, $i$ is a
kernel of $j$ and $j$ is a cokernel of $i$ in $\C$,
\item[b)] both classes $\mbold (\C)$ and $\ebold (\C)$ are subcategories
of $\C$; $\ebold (\C)$ is closed under base-changes along
arbitrary morphisms in $\C$ in the sense that for every exact
sequence $C' \to C \to C''$ and any morphism $f \colon D'' \to
C''$ in $\C$, there is a pullback commutative diagram
\[
\begin{CD}
C' @>>> D @>{j'}>> D''\\
@V{=}VV @VV{f'}V @VV{f}V \\
C' @>>> C @>{j}>> C''
\end{CD}
\]
where $j' \colon D \to D''$ is an admissible epimorphism; $\mbold
(\C)$ is closed under cobase-changes along arbitrary morphisms in
$\C$ in the (dual) sense that for every exact sequence $C' \to C
\to C''$ and any morphism $g \colon C' \to D'$ in $\C$, there is a
pushout diagram
\[
\begin{CD}
C' @>{i}>> C @>>> C'' \\
@V{g}VV @V{g'}VV @VV{=}V \\
D' @>{i'}>> D @>>> C''
\end{CD}
\]
where $i' \colon D' \to D$ is an admissible monomorphism,
\item[c)] if $f \colon C \to C''$ is a morphism with a kernel in $\C$,
and there is a morphism $D \to C$ so that the composition $D \to C
\to C''$ is an admissible epimorphism, then $f$ is an admissible
epimorphism; dually for admissible monomorphisms.
\end{enumerate}
According to Keller \cite{bK:90}, axiom (c) follows from the other
two.  We will use the standard notation $\rightarrowtail$ for
admissible monomorphisms and $\twoheadrightarrow$ for admissible
epimorphisms.
\end{DefRefName}

A \textit{preabelian category} is an additive category in which every morphism has a kernel and a cokernel \cite{aiG:96}.
Every morphism $f \colon F \to G$ in a preabelian category has a canonical decomposition
\[
f \colon X \, \xrightarrow{\ \coim (f) \ } \, \coim (f) \, \xrightarrow{ \ \bar{f} \ } \ \im (f) \, \xrightarrow{ \ \im (f) \ } \, Y
\]
where $\coim (f) =
\coker (\ker f)$ is the coimage of $f$ and $\im (f) = \ker (\coker f)$ is the image of $f$.
Recall that an \textit{abelian category} is a preabelian category
such that every morphism $f$ is
\textit{balanced}, that is, the canonical map
$\bar{f} \colon \coim (f) \to \im (f)$ is an
isomorphism.
An abelian category has the canonical exact structure
where all kernels and cokernels are respectively admissible
monomorphisms and admissible epimorphisms.

A \textit{subobject} of a fixed object $F$ is a monic $m \colon F'
\to F$.  The collection of all subobjects of $F$ forms a category
where morphisms are morphisms $j \colon F' \to F''$ between two
subobjects of $F$ such that $m'' j = m'$.  Notice that such $j$ are
also monic.  If the category is an exact category, there is the
subcategory of \textit{admissible subobjects} of $F$ represented by
admissible monomorphisms.  If both $m'$ and $m''$ are admissible, it
follows from exactness axiom 3 that $j$ is also an admissible
monomorphism.

Given two subobjects $m' \colon F' \to F$, $m'' \colon F'' \to F$,
the \textit{intersection} $F' \cap F''$, which is the pullback
of $m'$ along $m''$, is a subobject of $F$ and can be written as
the kernel of a morphism.  If $F'$ and $F''$ are
admissible subobjects then the intersection $F' \cap F''$ is
an admissible subobject.

\medskip

Now let $\A$ be a cocomplete abelian category.  The power set
$\mathcal{P}(X)$ of a proper metric space $X$ is partially ordered
by inclusion which makes it into the category with subsets of $X$ as
objects and unique morphisms $(S,T)$ when $S \subset T$. A
\define{presheaf} of $\A$-objects on $X$ is a functor $F \colon
\mathcal{P}(X) \to \A$. This corresponds to the usual notion of
presheaf of $\A$-objects on the discrete topological space
$X^\delta$ if the chosen Grothendieck topology on $\mathcal{P}(X)$
is the partial order given by inclusion, cf.\ section II.1 of
\cite{rH:77}. We will use terms which are standard in sheaf theory such as
\define{structure maps}, when referring to the morphisms $F(S,T)$.

\begin{DefRef}{FiltObjs}
A presheaf of objects in ${\bf A}$ on $X$ is an {\em $X$-filtered object
in ${\bf A}$} if all the structure maps of $F$ are monomorphisms.  We will
often suppress the reference to ${\bf A}$ when the meaning is clear from
the context.

For each presheaf $F$
there is an associated $X$-filtered object given by
\[
F_X (S) = \im F
(S,X).
\]

Suppose $F$ is an $X$-filtered object. Given a subobject $F'
\subset F(X)$ in $\A$, define the
\define{standard filtration} of $F'$ induced from $F$ by the
formula
\[
F'(S) = F(S) \cap F'.
\]
In other words, $F'(S)$ is the
image of the pullback
\[
\begin{CD}
 P @>>> F(S) \\
 @VVV @VVV \\
 F' @>>> F(X)
\end{CD}
\]
\end{DefRef}

\begin{DefRef}{BddCont}
The \textit{uncontrolled} category $\U (X,\A)$ is the category of
$X$-filtered objects in $\A$. The morphisms $F \to G$ in $\U (X,\A)$
are the morphisms $F(X) \to G(X)$ in $\A$.

Let $S[D]$ denote the subset $\{ x \in X \mid \dist(x,S) \le D \}$.
A morphism $f \colon F \to G$ in $\U (X, \A)$ is \textit{boundedly
controlled} if there is a number $D \ge 0$ such that the image of
$f$ restricted to $F(S)$ is a subobject of $G(S[D])$ for every
subset $S \subset X$.

The category $\Ub (X, \A)$ is
the full subcategory of $\U (X, \A)$ on the objects with the property $F (\varnothing) = 0$ and the
boundedly controlled morphisms.

If $f$ in addition has the property that for all subsets $S
\subset X$ the pullback $\im (f) \cap G(S)$ is a subobject of
$fF(S[D])$, then $f$ is called
\textit{boundedly bicontrolled}. In this case we say that $f$ has
filtration degree $D$ and write $\fil (f) \le D$.
\end{DefRef}

\begin{LemRef}{CatOK}
Let $f_1 \colon F \to G$, $f_2 \colon G \to H$ be in $\Ub (X, \A)$
and $f_3 = f_2  f_1$.
\begin{enumerate}
\item If $f_1$, $f_2$ are boundedly bicontrolled morphisms
and either $f_1 \colon F(X) \to G(X)$ is an epi or $f_2 \colon
G(X) \to H(X)$ is a monic, then $f_3$ is also boundedly
bicontrolled.
\item If $f_1$, $f_3$ are boundedly bicontrolled
and $f_1$ is epic then $f_2$ is also boundedly bicontrolled; if
$f_3$ is only boundedly controlled then $f_2$ is also boundedly
controlled.
\item If $f_2$, $f_3$ are boundedly bicontrolled
and $f_2$ is monic then $f_1$ is also boundedly bicontrolled; if
$f_3$ is only boundedly controlled then $f_1$ is also boundedly
controlled.
\end{enumerate}
\end{LemRef}

\begin{proof}
Suppose $\fil (f_i) \le D$ and $\fil (f_j) \le D'$ for $\{ i,j \}
\subset \{ 1, 2, 3 \}$, then in fact $\fil (f_{6-i-j}) \le D +
D'$ in each of the three cases. For example, there are
factorizations
\begin{gather}
f_2 G(S) \subset f_2 f_1 F(S[D]) =
f_3 F(S[D]) \subset H(S[D+D']) \notag \\
f_2 G(X) \cap H(S) \subset f_3 F(S[D']) = f_2 f_1 F(S[D'])
\subset f_2 G(S[D+D']) \notag
\end{gather}
which verify part 2 with $i=1$, $j=3$.
\end{proof}

\begin{PropRef}{BCAbelian}
$\Ub (X, \A)$ is an additive category with kernels and cokernels.
\end{PropRef}

\begin{proof}
Additive properties are inherited from $\A$. In particular, the
biproduct is given by the filtration-wise operation
\[
(F \oplus G)(S)
= F(S) \oplus G(S)
\]
in $\A$. For any boundedly controlled morphism
$f \colon F \to G$, the kernel of $f$ in $\A$ has the standard
$X$-filtration $K$ where
\[
K(S) = \ker (f) \cap F(S)
\]
which gives the
kernel of $f$ in $\Ub (X,\A)$. The canonical monic $\kappa \colon K
\to F$ has filtration degree $0$ and is therefore boundedly bicontrolled.
It follows from part 3 of Lemma
\refL{CatOK} that $K$ has the universal properties of the kernel in
$\Ub (X,\A)$.

Similarly, let $I$ be the standard $X$-filtration of the image of
$f$ in $\A$ by
\[
I(S) = \im (f) \cap G(S).
\]
Then there is a
presheaf $C$ over $X$ with
\[
C(S) = G(S)/I(S)
\]
for $S \subset X$.
Of course $C(X)$ is the cokernel of $f$ in $\A$. Recall fron Definition \refD{FiltObjs} that there
is an $X$-filtered object $C_X$ associated to $C$ given by
\[
C_X (S)
= \im C (S,X).
\]
The canonical morphism $\pi \colon G(X) \to C(X)$
gives a morphism of filtration $0$ (and which is therefore boundedly bicontrolled) $\pi \colon G \to C_X$ since
\[
\im (\pi  G(S,X)) = \im C(S,X) = C_X (S).
\]
This in conjunction with part 2 of Lemma \refL{CatOK} also verifies
the universal cokernel properties of $C_X$ and $\pi$ in $\Ub
(X,\A)$.
\end{proof}

\begin{RemRef}{NoAb}
If $\A$ is an abelian category and $X$ is unbounded then $\Ub (X,\A)$
is not necessarily an abelian category.

For an explicit description of a boundedly controlled morphism in
$\U (\mathbb{Z},\Mod (R))$ which is an isomorphism of left
$R$-modules but whose inverse is not boundedly controlled, we refer
to Example 1.5 of \cite{ePcW:85}.

This indicates that
under any embedding
of $\Ub (X,\A)$ in an abelian category $\F$ the kernels and
cokernels of some morphisms in $\F$ will be different from those in $\Ub
(X,\A)$.
\end{RemRef}

One consequence of Remark \refR{NoAb} is that $\Ub (X,\A)$ is not a
balanced category.

\begin{PropRef}{CharBB}
A morphism in $\Ub (X,\A)$ is balanced if and only
if it is boundedly bicontrolled.
\end{PropRef}

\begin{LemRef}{FutRef}
An isomorphism in $\Ub (X,\A)$ is boundedly bicontrolled.
\end{LemRef}

\begin{proof}
Suppose $f$ is an isomorphism in $\Ub (X,\A)$ bounded by $D(f)$ and let $f^{-1}$ be the inverse bounded by $D(f^{-1})$.  Then $f$ is boundedly bicontrolled with
filtration degree $\fil (f) \le \max \{ D(f^{-1}), D(f) \}$.
\end{proof}

\begin{LemRef}{mepseudo}
In any additive category, a morphism $h$ is monic if
and only if $\ker (h)$ exists and is the $0$ object. Similarly, $h$ is epic if
and only if $\coker (h)$ exists and is the $0$ object.
\end{LemRef}

\begin{proof}
Suppose $h_1$, $h_2 \colon F \to G$ and $h \colon G \to H$ are
such that $h h_1 = h h_2$, then $h (h_1 - h_2) = 0$. So there is
a morphism $F \to \ker (h) = 0$ such that $F \to \ker (h) \to G$
is precisely $h_1 - h_2$. Hence $h_1 - h_2 = 0$ and $h_1 = h_2$.
Conversely, if $h$ is monic in a category with a zero object, it is
clear that $\ker (h) = 0$. The fact about epics is similar.
\end{proof}

\begin{proof}[Proof of Proposition \refP{CharBB}]
Let $f \colon F \to G$ be a morphism in $\Ub (X,\A)$, and let $J$ be the coimage and $I$ be the image of $f$.
The standard filtration $I(S) = I \cap G(S)$ makes the inclusion $i \colon I \to G$ boundedly bicontrolled of filtration $0$.  Similarly, $J$ is the cokernel of the inclusion of $\ker (f)$ in $F$, so $J$ has the filtration described in the proof of Proposition \refP{BCAbelian} which makes the projection $p \colon F \to J$ boundedly bicontrolled of filtration $0$.

Necessity of the condition follows from Lemma \refL{FutRef}.

Now $f$ factors as the composition
\[
F \ \xrightarrow{\ p \ } \, J \, \xrightarrow{\ \theta \ } \, I \, \xrightarrow{\ i \ } \, G,
\]
where $\theta$ is the canonical map.  If $f$ is bounded by $D$ then clearly $\theta$ is bounded by $D$
and has the $0$ object for the kernel and the cokernel.  This shows that $\theta$ is an isomorphism in $\A$
by Lemma \refL{mepseudo}.  In particular, there is an inverse $\theta^{-1} \colon J \to I$
in $\A$.  Now the condition that
\[
I \cap G(S) \subset fF(S[b])
\]
for some number $b$ and a subset $S \subset X$ is equivalent to the condition
\[
\theta^{-1} (I(S)) \subset J(S[b]).
\]
So $f$ is
boundedly bicontrolled if and only if $\theta^{-1}$ is bounded and, therefore, is an isomorphism in $\Ub (X,\A)$.
\end{proof}

\begin{CorRef}{IsoBB}
An isomorphism in $\Ub (X,\A)$ is a morphism which is an isomorphism in $\A$ and is boundedly bicontrolled.
\end{CorRef}

\begin{DefRef}{ExStrUb}
The \textit{admissible monomorphisms} $\mUb (X,\A)$ in $\Ub (X,\A)$
consist of boundedly bicontrolled morphisms $m \colon F_1 \to F_2$
such that $m \colon F_1 (X) \to F_2 (X)$ is a monic in $\A$. The
\textit{admissible epimorphisms} $\eUb (X,\A)$ are the boundedly
bicontrolled morphisms $e \colon F_1 \to F_2$ such that $e \colon
F_1 (X) \to F_2 (X)$ is an epi in $\A$. The class $\mathcal{E}$ of
exact sequences consists of the sequences
\[
E^{\subdot} \colon \quad E' \xrightarrow{\ i \ } E \xrightarrow{\ j \ } E''
\]
with $i \in \mUb (X,\A)$ and $j \in \eUb (X,\A)$ which are exact at
$E$ in the sense that $\im (i)$ and $\ker (j)$ represent the same
subobject of $E$.
\end{DefRef}

\begin{ThmRef}{UbisWAb}
$\Ub (X,\A)$ is an exact category.
\end{ThmRef}

\begin{proof}
(a) It follows from Lemma \refL{CatOK} that any short exact
sequence $F^{\subdot}$ isomorphic to some $E^{\subdot} \in
\mathcal{E}$ is also in $\mathcal{E}$, that
\[
\begin{CD}
F' @>{[\id,0]}>> F' \oplus F'' @>{[0,\id]^T}>> F''
\end{CD}
\]
is in $\mathcal{E}$, and that $i = \ker (j)$, $j = \coker (i)$ in
any $E^{\subdot} \in \mathcal{E}$.

(b) The collections of morphisms $\mUb (X,\A)$ and $\eUb (X,\A)$ are
closed under composition by part 1 of Lemma \refL{CatOK}.

Now suppose we are given
an exact sequence
\[
E^{\subdot} \colon \quad E' \xrightarrow{\ i \ } E \xrightarrow{\ j \ } E''
\]
in $\mathcal{E}$ and a morphism $f \colon A \to E'' \in \Ub
(X,\A)$.  Let $D(j)$ be a filtration constant for $j$ as a boundedly controlled epi and let $D(f)$ be a bound for $f$ as a boundedly controlled map.
There is a base change diagram
\[
\begin{CD}
E' @>>> E \times_{f} A @>{j'}>> A\\
@V{=}VV @VV{f'}V @VV{f}V \\
E' @>>> E @>{j}>> E''
\end{CD}
\]
where $m \colon E \times_{f} A \to E \oplus A$ is the kernel of
the epi
\[
j  \circ \pr_1 - f  \circ \pr_2 \colon E \oplus A \longrightarrow E''
\]
and $f' =
\pr_1  \circ \, m$, $j' = \pr_2  \circ \, m$.
The $X$-filtration on $E \times_{f} A$ is the standard filtration
\[
\left( E \times_{f} A \right) (S) = E \times_{f} A \cap \left(
E(S) \times A(S) \right).
\]
The induced map $j'$ has the same kernel as
$j$ and is bounded by $0$ since
\[
j' \left( (E \times_{f} A) (S) \right) \subset A(S).
\]
In fact,
\[
fA(S) \subset E''(S[D(f)]),
\]
so
\[
fA(S) \subset j E(S[D(f) + D(j)]),
\]
and
\[
\im (j') \cap A(S) \subset j' \left( E \times_{f} A \right) \left(
S[D(f) + D(j)] \right).
\]
This shows that $j'$ is boundedly
bicontrolled of filtration degree $D(f) + D(j)$.

Therefore, the class of admissible epimorphisms is
closed under base change by arbitrary morphisms in $\Ub (X,\A)$.
Cobase changes by admissible monomorphisms are similar.
\end{proof}

The following Proposition is an organic characterization of the exact structure in $\Ub (X,\A)$.

\begin{PropRef}{SSSST}
The exact structure $\mathcal{E}$ in $\Ub (X,\A)$ consists of sequences isomorphic to those
\[
E^{\subdot} \colon \quad E' \xrightarrow{\ i \ } E \xrightarrow{\ j \ } E''
\]
which possess restrictions
\[
E^{\subdot} (S) \colon \quad E' (S) \xrightarrow{\ i \ } E (S) \xrightarrow{\ j \ } E'' (S)
\]
for all subsets $S \subset X$, and each $E^{\subdot} (S)$ is an exact sequence in $\A$.
\end{PropRef}

\begin{proof}
Clearly, each of the sequences $E^{\subdot}$ described in the statement is an exact sequence in $\Ub (X,\A)$.  Indeed, the restriction $i \colon E'(S) \to E(S)$ is monic and $j \colon E(S) \to E''(S)$ is epic, so $i \colon E' \to E$ and $j \colon E \to E''$ are both bicontrolled of filtration $0$.

Suppose $F^{\subdot}$ is a sequence isomorphic to such $E^{\subdot}$.  There is a commutative diagram
\[
\begin{CD}
F' @>{f}>> F @>{g}>> F''\\
@V{\cong}VV @VV{\cong}V @VV{\cong}V \\
E' @>{i}>> E @>{j}>> E''
\end{CD}
\]
Then $f$ and $g$ are compositions of two isomorphisms (which are boundedly bicontrolled by Lemma \refL{FutRef}) which are either preceded by a boundedly bicontrolled monic or followed by a boundedly bicontrolled epi.  By Lemma \refL{CatOK}, part (1), both  $f$ and $g$ are boundedly bicontrolled.

Now suppose $F^{\subdot}$ is an exact sequence in $\mathcal{E}$.  Let $K = \ker (g)$ and $C = \coker (f)$, then we obtain a commutative diagram
\[
\begin{CD}
F' @>{f}>> F @>{g}>> F''\\
@V{\cong}VV @VV{=}V @AA{\cong}A \\
K @>{i}>> F @>{j}>> C
\end{CD}
\]
where the vertical maps are the canonical isomorphisms.
By the construction of kernels and cokernels in Proposition \refP{BCAbelian}, there are exact sequences
\[
K(S) \xrightarrow{\ i \ } F(S) \xrightarrow{\ j \ } C(S)
\]
for all subsets $S \subset X$.
\end{proof}

\begin{DefRef}{ClUnExt}
A full subcategory $\mathbf{H}$ of a small exact category $\C$ is
said to be \textit{closed under extensions} in $\C$ if $\mathbf{H}$
contains a zero object and for any exact sequence $C' \to C \to C''$
in $\C$, if $C'$ and $C''$ are isomorphic to objects from
$\mathbf{H}$ then so is $C$.

A \textit{Grothendieck subcategory} of
an exact category is a subcategory which is closed under
isomorphisms, exact extensions, admissible subobjects, and
admissible quotients.
\end{DefRef}

It is known \cite{dQ:73} that a subcategory closed under
extensions inherits an exact structure from $\C$.

Now let $\E$ be a Grothendieck subcategory of $\A$ and let $F$ be an
object of $\Ub (X,\A)$.

\begin{DefRef}{RealBCprelim}
(1) $F$ is $\E$-\textit{local} if $F (V)$ is an object of $\E$ for
every bounded subset $V \subset X$.

(2) $F$ is \textit{lean} or $D$-\textit{lean} if there is a number
$D \ge 0$ such that for every subset $S$ of $X$
\[
F(S) \subset \sum_{x \in S} F(B_D (x)),
\]
where $B_D (x)$ is the metric ball of radius $D$ centered at $x$.

(3) $F$ is \textit{insular} or $d$-\textit{insular} if there is a
number $d \ge 0$ such that
\[
F(T) \cap F(U) \subset F(T[d] \cap U[d])
\]
for every pair of subsets $T$, $U$ of $X$.

Notice that a $d$-insular object has the property that for any
subset $T \subset X$,
\[
F(T) \cap F(U) = 0
\]
whenever $T \cap U[2d]
= \varnothing$.
\end{DefRef}

\begin{RemRef}{leaninscld}
It is clear that properties (1), (2), and (3) are preserved under
isomorphisms in $\Ub (X,\A)$.
\end{RemRef}

\begin{PropRef}{lninpres}
$\mathrm{(1)}$  Lean objects are closed under exact extensions in
$\Ub (X,\A)$, that is, if
\[
E' \longrightarrow E \longrightarrow E''
\]
is an exact sequence in $\Ub (X,\A)$, and $E'$, $E''$ are lean, then
$E$ is lean.

$\mathrm{(2)}$ Insular objects are closed under exact extensions in
$\Ub (X,\A)$.

$\mathrm{(3)}$ If in the exact sequence above the object $E$ is
lean and insular then

$\quad \mathrm{(a)}$ $E'$ is insular,

$\quad \mathrm{(b)}$ $E''$ is lean,

$\quad \mathrm{(c)}$ $E''$ is insular if and only if $E'$ is lean.
\end{PropRef}

\begin{proof}
Let
\[
E' \xrightarrow{\ f \ } E \xrightarrow{\ g \ } E''
\]
be an exact sequence in $\Ub (X, \A)$ and let $b \ge 0$ be a common
filtration degree for both $f$ and $g$ as boundedly bicontrolled
maps.

(1) Assume that both $E'$ and $E''$ are $D$-lean as objects of $\Ub
(X, \A)$. Consider $E(S)$, then
\[
gE(S) \subset E'' (S[b])
\]
and so
\[
gE(S) \subset \sum_{x \in S[b]} E'' (B_D (x)).
\]
For each $x \in S[b]$,
\[
E'' (B_D (x)) \subset gE(B_{D + b} (x)),
\]
so
\[
E(S) \subset \sum_{x \in S[b]} E (B_{D+2b} (x)) + \sum_{x \in S[b]}
fE' (B_{D + 2b} (x)).
\]
Therefore
\[
E(S) \subset \sum_{x \in S[b]} E (B_{D + 3b} (x)) \subset \sum_{x
\in S} E (B_{D + 4b} (x)),
\]
so $E$ is $(D + 4b)$-lean.

(2) Assuming that both $E'$ and $E''$ are $d$-insular, for any pair
of subsets $T$ and $U$ of $X$,
\begin{equation} \begin{split}
&g (E(T) \cap E(U))\\
\subset\ &E'' (T[b]) \cap E'' (U[b])\\
\subset\ &E'' (T[b + d] \cap U[b + d]).
\end{split} \notag \end{equation}
Now we have
\begin{equation} \begin{split}
&E(T) \cap E(U)\\
\subset\ &E(T[2 b + d] \cap U[2 b + d]) + fE' \cap E(T[2b +d]) \cap E (U[2b + d])\\
\subset\ &E(T[2 b + d] \cap U[2 b + d]) + f (E' (T[3 b + d]) \cap E' (U[3 b + d])\\
\subset\ &E(T[2 b + d] \cap U[2 b + d]) + f E' (T[3 b + 2d] \cap U[3 b + 2d])\\
\subset\ &E(T[4 b + 2d] \cap U[4b + 2d]).
\end{split} \notag \end{equation}
So $E$ is $(4b + 2d)$-insular.

(3a) Suppose $E$ is $d$-insular.  Given subsets $T$ and $U$ of
$X$,
\begin{equation} \begin{split}
&f (E'(T) \cap E'(U))\\
\subset\ &fE'(T) \cap fE'(U)\\
\subset\ &E(T[b]) \cap E(U[b])\\
\subset\ &E(T[b + d]) \cap E(U[b + d]),
\end{split} \notag \end{equation}
so
\[
E'(T) \cap E'(U) \subset E' (T[2b+d] \cap U[2b+d]).
\]
Thus $E'$ is $(2b + d)$-insular.

(3b) If $E$ is $D$-lean then for any $S \subset X$, $E'' (S)
\subset gE(S[b])$.  Since
\[
E(S[b]) \subset \sum_{x \in S[b]} E (B_{D} (x)),
\]
\[
E'' (S) \subset \sum_{x \in S[b]} E'' (B_{D + b} (x)).
\]
So
\[
E'' (S) \subset \sum_{x \in S} E'' (B_{D + 2b} (x)),
\]
and $E''$ is $(D+2b)$-lean.

(3c) Suppose $E'$ is $D$-lean and $E$ is $d$-insular.  For any pair of subsets $T$, $U
\subset X$,
\[
E'' (T) \cap E'' (U) \subset gE(T[b]) \cap gE(U[b]).
\]
Given
\[
z \in E'' (T) \cap E'' (U),
\]
let $y' \in E(T[b])$ and $y''
\in E(U[b])$ so that
\[
g(y') = g(y'') = z.
\]
Now
\[
k = y' - y'' \in \big( E(T[b]) + E(U[b]) \big) \cap \ker (g),
\]
so there is
\[
\overline{k} \in E' (T[2b]) + E'(U[2b]) \subset E' (T[2b] \cup U[2b])
\]
with $f(\overline{k}) =k$.
Since $E'$ is $D$-lean,
\[
\overline{k} \in \sum_{x \in T[2b] \cup U[2b]} E' (B_D (x)) =
\sum_{x \in T[2b]} E' (B_D (x)) + \sum_{y \in U[2b]} E' (B_D (y)).
\]
Hence,
\[
\overline{k} \in E' (T[2b+D]) + E' (U[2b+D]).
\]
This allows us to write $\overline{k} = \overline{k}_1 +
\overline{k}_2$, where $\overline{k}_1 \in E' (T[2b+D])$ and
$\overline{k}_2 \in E' (U[2b+D])$.  Now
\[
k = f \overline{k}_1 + f
\overline{k}_2.
\]
Notice that
\[
y' = y'' + k = y'' + f \overline{k}_1 + f \overline{k}_2.
\]
So
\[
y = y' - f \overline{k}_1 = y'' + f \overline{k}_2
\]
has the property
\[
y \in E (T [3b+D]) \cap E(U[3b+D]) \subset E (T [3b+D+d] \cap
U[3b+D+d]),
\]
and $g(y) = z$. Hence
\[
z \in E'' (T[4b+D+d] \cap U[4b+D+d]).
\]
We conclude that $E''$ is $(4b+D+d)$-insular. The converse is proved
similarly; it is not used in this paper.
\end{proof}

\begin{DefRef}{strict}
An object $F$ of $\Ub (X,\A)$ is called $\ell$-\textit{strict} or simply \textit{strict} if there exists an order preserving function
\[
\ell \colon \mathcal{P}(X) \longrightarrow [0, \infty)
\]
such that for every subset $S$ of $X$ the subobject
$F_{S} = F(S)$
is
$\E$-local, $\ell_{S}$-lean and $\ell_{S}$-insular with respect to the
standard filtration
$F_{S} (U) = F_{S} \cap F(U)$.
\end{DefRef}

Unlike the subcategory of lean and insular objects, the subcategory of strict
objects is not necessarily closed under isomorphisms.

\begin{DefRef}{RealBC}
The \define{boundedly controlled} category $\B (X,\E)$ is the full
subcategory of $\Ub (X,\A)$ on objects which are isomorphic to
strict objects.
\end{DefRef}

The terminology adopted here is convenient and should not suggest
relations to boundedly controlled spaces and maps introduced earlier
by Anderson and Munkholm \cite{dAhM:88}.

\begin{RemRef}{CocNot}
The exact subcategory $\E$ is not assumed to be cocomplete. In fact,
the construction is most interesting when it is not. Notice also
that the notation $\B (X, \E)$ does not suggest that the objects $F$
have the terminal piece $F(X)$ in $\E$, unlike the notation for $\Ub
(X, \A)$ where $F(X)$ are in $\A$.  The object $F(X)$ is contained
in the cocompletion of $\E$ in $\A$.
\end{RemRef}

\begin{ThmRef}{ExtCl}
$\B (X,\E)$ is closed under extensions in $\Ub (X,\A)$.
\end{ThmRef}

\begin{proof}
Let
\[
E' \xrightarrow{\ f \ } E \xrightarrow{\ g \ } E''
\]
be an exact sequence in $\Ub (X, \E)$ and let $b \ge 0$ be a common
filtration degree for both $f$ and $g$ as boundedly bicontrolled
maps.  We will also assume, without loss of generality, that both
$E'$ and $E''$ are $\ell$-strict for some function $\ell \ge 0$. We
need to check that $E$ is isomorphic to a strict object.

Since $\E$ is a Grothendieck subcategory of $\A$, for every bounded
subset $V \subset X$ the restriction
\[
g \vert E(V) \colon E(V) \longrightarrow
gE(V)
\]
is an admissible epimorphism onto an admissible subobject of
$E''(V[D])$, which is in $\E$. The kernel of $g \vert E(V)$ is the
admissible subobject $\ker (g) \cap E(V)$ of $E(V)$, which is also
in $\E$.  So $E(V)$ is in $\E$ by closure under extensions in $\A$.

To see that $E$ is isomorphic to a strict object,
consider $S \subset X$ so that ${E}''(S[b])$
is $\ell_{S[b]}$-lean and $\ell_{S[b]}$-insular. The induced epi
\[
g \colon E(S[2b]) \cap g^{-1} {E}''(S[b]) \longrightarrow
{E}''(S[b])
\]
extends to another epi
\[
g' \colon f {E}'(S[3b]) + E(S[2b]) \cap g^{-1} {E}''(S[b])
\longrightarrow {E}''(S[b])
\]
with $\ker (g') = {E}'(S[3b])$.

Since both ${E}'(S[3b]$ and ${E}''(S[b])$
are $\E$-local, $\ell_{S[3b]}$-lean, and $\ell_{S[3b]}$-insular,
parts (1) and (2) of Proposition
\refP{lninpres} show that
\[
\widehat{E}(S) = f {E}'(S[3b]) + E(S[2b]) \cap g^{-1} {E}''(S[b])
\]
is $(3b + 4b)$-lean and $(6b+4b)$-insular.
This makes
the filtration $\widehat{E}$ $\phi$-strict for the function $\phi_S = \ell_{S[10b]}$.

Clearly, the identity map $\id \colon E(X) \to E(X)$ gives an
isomorphism $\id \colon \widehat{E} \to E$ with $\fil (\id) \le 4b$.
\end{proof}

\begin{CorRef}{BisWAb}
$\B (X,\E)$ is an exact category in the sense of Quillen.
\end{CorRef}

The bounded theory of geometric free modules described in the Introduction can be generalized to arbitrary additive categories.

\begin{DefRef}{UIYWA}
Given a proper metric space $M$ and an additive category $\mathcal{A}$,
Pedersen--Weibel \cite{ePcW:85} define the category of \textit{geometric objects}
$\mathcal{B} (M, \mathcal{A})$ as follows.
The objects are functions $F$ from $M$ to the objects of $\mathcal{A}$
which satisfy the local finiteness condition: a bounded subset of $M$ contains only finitely many points $x \in M$ such that the values $F(x)$ are nonzero.
A morphism $\phi \colon F \to G$ of degree $D \ge 0$ is a collection of $\mathcal{A}$-morphisms
\[
\phi (x,y) \colon F(x) \to G(y)
\]
such that $\phi (x,y)$ is zero unless $d(x,y) \le D$.

The category $\mathcal{B} (M, \mathcal{A})$ is an additive category with the biproduct
\[
(F \oplus G)(x) = F(x) \oplus G(x).
\]
The associated bounded $K$-theory spectrum is denoted by $K (M, \mathcal{A})$.
\end{DefRef}

Given an exact category $\E$, one can view $\E$ as an additive category with the underlying split exact structure.  When we use $\E$ as coefficients in a category of geometric objects $\mathcal{B} (M, \E)$,
this is the structure that is implicitly understood.

\begin{CorRef}{BisWAb2}
The
additive category $\mathcal{B}(X,\E)$ of geometric objects with
the split exact structure is an exact subcategory of $\B (X,\E)$.
\end{CorRef}

\begin{proof}
The $X$-filtration of the geometric objects in $\mathcal{B}(X,\E)$
is given by
\[
F(S) = \bigoplus_{x \in S} F(x),
\]
and the
structure maps are the boundedly controlled inclusions and
projections onto direct summands.
\end{proof}

Recall that a morphism $e \colon F \to F$ is an idempotent if
$e^2=e$. Categories in which every idempotent is the projection onto
a direct summand of $F$ are called \textit{idempotent complete}.
Abelian categories are clearly idempotent complete. Thus $\A$ and
its Grothendieck subcategories, which are abelian, are idempotent
complete.

\begin{CorRef}{IdCIdC}
The subcategory $\B (X,\E)$ is idempotent complete.
\end{CorRef}

\begin{proof}
Since the restriction of an idempotent $e$ to the image of $e$ is
the identity, every idempotent is boundedly bicontrolled of
filtration $0$. It follows easily that the splitting of $e$ in $\A$
is in fact a splitting in $\B (X,\E)$.
\end{proof}

\begin{PropRef}{ClUCok}
The subcategory $\B (X,\E)$ is closed under admissible quotients of strict objects.  Precisely, for
a given boundedly bicontrolled epi $f \colon F \to G$ in $\Ub (X,\A)$ where both $F$ and the kernel $k \colon K \to F$ with the standard filtration $K(S) = K \cap F(S)$ are strict, the cokernel $G$ is isomorphic to a strict object.
\end{PropRef}

\begin{proof}
Suppose $\fil (f) \le b$, then from the assumptions
\[
K(S[b]) \longrightarrow F(S[b]) \longrightarrow fF(S[b])
\]
is an exact sequence in $\Ub (X,\A)$ for any subset $S \subset X$.
Since $F(S[b])$ is lean and insular and $K(S[b])$ is lean, the quotient $fF(S[b])$ is lean and insular by Proposition \refP{lninpres}.  It is clear that $fF(S[b])$ is also $\E$-local.  Thus the object $\widehat{G}$ with filtration
\[
\widehat{G}(S)=fF(S[b])
\]
is strict.  The identity map induces an isomorphism $\id \colon G \to
\widehat{G}$ with $\fil (\id) \le 2b$ because for all $S \subset X$ we have $G(S) \subset
\widehat{G}(S)$ and $\widehat{G}(S) \subset G(S[2b])$.
\end{proof}

\begin{RemRef}{NHYW}
For additional flexibility, one may want to impose weaker
requirements on objects in $\Ub (X,\A)$.  Restricting as in
Definition \refD{RealBC} to objects $F$ with a fixed locally finite
covering $\mathcal{U} \subset \mathcal{P}(X)$ by bounded subsets $U
\in \mathcal{B}(X)$ such that
\[
F(X) = \sum_{U \in \mathcal{U}} F(U)
\]
gives another exact category. In this case, one may also relax the
bounded control conditions on the maps to those of Lipschitz type.
Similar modifications have become useful in recent work of
Hambleton--Pedersen \cite{iHeP:99} and Pedersen--Weibel
\cite{ePcW:97} in controlled $K$-theory.
\end{RemRef}

\SecRef{Localization in controlled categories}{QCC}

\begin{DefRef}{Germs}
Let $F$ be an object of $\B (X, \E)$ and $Z$ be a subset of $X$.
We say $F$ is \textit{supported near} $Z$ if there is a number $d
\ge 0$ such that $F(X) \subset F(Z[d])$.

Let $\B (X,\E)_{<Z}$ be the full subcategory of $\B (X,\E)$ on
objects supported near $Z$. If $\B_{d} (X,\E)_{<Z}$ denotes the full
subcategory of $\B (X,\E)$ with objects $F$ as above then
\[
\B (X,\E)_{<Z} = \colim{d} \B_{d} (X,\E)_{<Z}.
\]
\end{DefRef}

\begin{PropRef}{Serre}
$\B (X,\E)_{<Z}$ is a Grothendieck subcategory of $\B (X,\E)$.
\end{PropRef}

\begin{proof}
First we show closure under exact extensions. Let
\[
F' \xrightarrow{\ f\ } F \xrightarrow{\ g\ } F''
\]
be an exact sequence in $\B (X,\E)$. Let $b$ be a common filtration
degree of $f$ and $g$
and let $d'$, $d'' \ge 0$ be numbers with $F' = F' (Z[d'])$ and $F''
= F'' (Z[d''])$. Since $F = I + M$, where $I = \im (f)$ and $M$ is
any subobject $M \subset F$ with $g (M) = F''$, it suffices to show
that for some $d \ge 0$
\[
I (X) = I (Z[d]) \subset F(Z[d]),
\]
and that $M$ can be chosen to be a subobject of $F(Z[d])$. Indeed,
\begin{gather}
I (X) = f F'(X) = f F'(Z[d']) \subset F(Z[d'+b]),
\notag \\
F'' (X) = g F(X) \cap F'' (Z[d'']) \subset g F(Z[d''+b]). \notag
\end{gather}
Let $M = F (Z[d''+b])$.  If we choose $d = \max \{ d'+b, d''+b \}$
then $F = F(Z[d])$ is in $\B (X,R)_{<Z}$.

Now suppose $f \colon F' \to F$ is an admissible subobject in $\B
(X,\E)$, which is a boundedly bicontrolled monic with $\fil (f) \le
b$, $F = F(Z[d])$, and $F$ is $c$-insular. Also suppose $F$ and $F'$
are respectively $D$- and $D'$-lean, then notice that
\[
f F' (B_{D'}
(x)) \subset F(B_{D' + b} (x)),
\]
while $F(X) \subset F(Z[d])$.
Therefore,
\[
F' (B_{D'} (x)) = 0
\]
for
\[
x \in X - Z[d + D + D' + b +
2c].
\]
This means that
\[
F' = F'(Z[d + D + 2D' + b + 2c]).
\]
On the
other hand, if $g \colon F \to F''$ is an admissible quotient of
filtration $b$ then it is easy to see that
\[
F'' = F''(Z[d + D + b])
\]
is also in $\B (X,\E)_{<Z}$.  Since $\B (X,\E)_{<Z}$ is clearly
closed under isomorphisms, this proves the assertion.
\end{proof}

Given an object $F \in \B (X, \E)$ and a subset $T \subset X$,
we will need a construction of an admissible subobject
$\widetilde{F}$ of $F$ in $\B (X, \E)$ such that
\[
F(T) \subset \widetilde{F} \subset F(T[D])
\]
for some $D \ge 0$.

Choose a strict $F'$ isomorphic to $F$ in $\B (X, \E)$ and assume the chosen
isomorphism and its inverse are of filtration $b$.

\begin{LemRef}{Help}
The subobject $\widetilde{F} = F'(T[b])$ is an admissible
subobject of $F$ in $\B (X, \E)$ and satisfies
\[
F(T) \subset \widetilde{F} \subset F(T[2b]).
\]
\end{LemRef}

\begin{proof}
The cokernel $G'$ of the inclusion $k \colon F'(T[b]) \to F$ is in $\B (X,\E)$ by Proposition
\refP{ClUCok}.  We can view $F'(T[b])$ as an admissible subobject of $F$ with the cokernel $G$ isomorphic to $G'$.
\end{proof}

\begin{DefRef}{CLFracs}
A class of morphisms $\Sigma$ in an additive category $\C$
\textit{admits a calculus of right fractions} if
\begin{enumerate}
\item the identity of each object is in $\Sigma$,
\item $\Sigma$ is closed under composition,
\item each diagram $F \xrightarrow{\ f} G \xleftarrow{\ s\ } G'$ with $s
\in \Sigma$ can be completed to a commutative square
\[
\begin{CD}
 F' @>{f'}>> G'\\
 @VV{t}V @VV{s}V\\
 F @>{f}>> G
\end{CD}
\]
with $t \in \Sigma$, and
\item if $f$ is a morphism in $\C$ and $s \in \Sigma$ such that
$sf = 0$ then there exists $t \in \Sigma$ such that $ft = 0$.
\end{enumerate}
In this case there is a construction of the \textit{localization}
$\C [\Sigma^{-1}]$ which has the same objects as $\C$.  The
morphism sets $\Hom (F,G)$ in $\C [\Sigma^{-1}]$ consist of
equivalence classes of diagrams
\[
(s,f) \colon \quad F \xleftarrow{\ s\ } F' \xrightarrow{\ f} G
\]
with the equivalence relation generated by $(s_1,f_1) \sim
(s_2,f_2)$ if there is a map $h \colon F'_1 \to F'_2$ so that $f_1
= f_2 h$ and $s_1 = s_2 h$.
Let $(s \vert f)$ denote the
equivalence class of $(s,f)$.

The composition of morphisms in $\C
[\Sigma^{-1}]$ is defined by
\[
(s \vert f) \circ (t \vert g) = (st' \vert gf')
\]
where $g'$ and $s'$ fit in the commutative square
\[
\begin{CD}
 F'' @>{f'}>> G'\\
 @VV{t'}V @VV{t}V\\
 F @>{f}>> G
\end{CD}
\]
from axiom 3.
\end{DefRef}

\begin{PropRef}{FactsFrac}
The localization $\C [\Sigma^{-1}]$
is a category.  The morphisms
of the form $(\id \vert s)$ where
$s \in \Sigma$ are isomorphisms
in $\C [\Sigma^{-1}]$. The rule
$P_{\Sigma} (f) = (\id \vert f)$
gives a functor
\[
P_{\Sigma} \colon \C \longrightarrow
\C [\Sigma^{-1}]
\]
which
is universal among the functors making
the morphisms $\Sigma$
invertible.
\end{PropRef}

\begin{proof}
The proofs of these facts can be found in Chapter I of
\cite{pGmZ:67}.  The inverse of $(\id \vert s)$ is $(s \vert
\id)$.
\end{proof}

Suppose $\E$ is a Grothendieck subcategory of a cocomplete abelian
category $\A$, and let $\Z$ be the subcategory $\B (X, \E)_{<Z}$ of
$\B = \B (X, \E)$ for a fixed choice of $Z \subset X$. Let the class
of \textit{weak equivalences} $\Sigma$ consist of all finite
compositions of admissible monomorphisms with cokernels in $\Z$ and
admissible epimorphisms with kernels in $\Z$. We will show that the
class $\Sigma$ admits a calculus of right fractions.

\begin{DefRef}{rfil}
A Grothendieck subcategory $\Z \subset \B$ is \textit{right
filtering} if each morphism $f \colon F \to G$ in $\B$, where $G$
is an object of $\Z$, factors through an admissible epimorphism $e
\colon F \to \overline{G}$, where $\overline{G}$ is in $\Z$.
\end{DefRef}

\begin{LemRef}{SubSubLem}
The subcategory $\Z = \B(X,\E)_{<Z}$ of $\B = \B(X,\E)$ is right filtering.
\end{LemRef}

\begin{proof}
Suppose first that $F$ and $G$ are strict with the characteristic functions
$\ell _F$ and $\ell _G$ respectively.
Let $L_F = \ell _F (X)$ and $L_G = \ell _G (X)$.
If
\[
G = G(Z[d_G])
\]
and the given morphism $f \colon F \to G$ is bounded by $d$
then we have
\[
f F(B_{L_G}(x)) \subset
G(B_{{L_G} + b} (x)) = 0
\]
for all
\[
x \in X - Z[d_G + L_G + 2L_G + d + L_F].
\]
Let
\[
E = F(X - Z[d_G + L_G + 2L_G + d + L_F]),
\]
then $f E =0$.
Now $E$ is an admissible subobject of $F$ by Lemma \refL{Help};
let $\overline{G}$ be the
cokernel of the inclusion.  Since
\[
\overline{G} (B_{L_F}(x)) =0
\]
for all
\[
x \in X
- Z[d_G + L_G + 2L_G + d + L_F],
\]
we have
\[
\overline{G} = \overline{G} (Z[d_G + 2L_G + 2L_G + d + L_F])
\]
as an object of $\Z$.
The required factorization is the right square in the map between the two exact sequences
\[
\begin{CD}
E @>>> F @>{j'}>> \overline{G}\\
@V{i}VV @VV{=}V @VVV \\
K @>{k}>> F @>{f}>> G
\end{CD}
\]

If $F$ and $G$ are not strict, one considers a map $f' \colon F' \to G'$ between strict objects isomorphic to $F$ and $G$ and chooses
the subobject
\[
E = F'(X - Z[d_G + L_G + 2L_G + d + L_F + 4b])
 \]
 of $F'$ for an appropriate value of $b$, as in the proof of Lemma \refL{Help}.
\end{proof}

\begin{CorRef}{WeakEqs}
The class $\Sigma$ admits a calculus
of right fractions.
\end{CorRef}

\begin{proof}
This follows from Lemma \refL{SubSubLem}, see Lemma 1.13 of \cite{mS:03}.
\end{proof}

\begin{DefRef}{Quot}
The
\define{quotient category}
$\mathbf{B}/\mathbf{Z}$ is the localization $\B [\Sigma^{-1}]$.
\end{DefRef}

It is clear that the quotient $\mathbf{B}/\mathbf{Z}$ is an
additive category, and $P_{\Sigma}$ is an additive functor.
In fact, we have the following.

\begin{ThmRef}{ExLocStr}
The short sequences in $\mathbf{B}/\mathbf{Z}$ which are
isomorphic to images of exact sequences from $\B$ form a Quillen
exact structure.
\end{ThmRef}

The proof uses the following fact.

\begin{DefRef}{Schkiy}
An extension closed full subcategory $\Z$ of an exact category $\B$ is called \textit{right s-filtering}
if whenever $f \colon F \to G$ is an admissible
monomorphism with $F$ in $\Z$ then there exist an object $E$ in $\Z$ and an
admissible epimorphism $e \colon G \to E$ such that the composition $ef$ is an
admissible monomorphism.
\end{DefRef}

\begin{LemRefName}{KOIUYTR}{Schlichting}
Let $\Z$ be a right filtering and right s-filtering subcategory in $\B$.  Then the quotient category $\mathbf{B}/\mathbf{Z}$, equipped with the images of the exact sequences from $\B$, is an exact category.  Moreover, exact functors from $\B$ vanishing on $\Z$ are in bijective correspondence with exact functors from $\mathbf{B}/\mathbf{Z}$.
\end{LemRefName}

\begin{proof}
See Proposition 1.16 of \cite{mS:03}.
\end{proof}

\begin{proof}[Proof of Theorem \refT{ExLocStr}]
Since $\Z$ is right filtering by Lemma \refL{SubSubLem}, it remains
to check that $\Z$ is right s-filtering.

Again, suppose that $F$ and $G$ are strict with the characteristic functions
$\ell _F$ and $\ell _G$ and let $L_F = \ell _F (X)$ and $L_G = \ell _G (X)$.  Assume that $F = F(Z[d_F])$, $\fil (f) \le d$, and let
\[
G' = G(X-Z[d_F + 2L_F + 2L_G + d + L_G]).
\]
Let $e \colon G \to E$ be the cokernel of the inclusion, then $f
(F) \cap G' = 0$, so that $ef$ is an admissible
monomorphism with $\fil (ef) = \fil (f) \le d$.  If $F$ and $G$ are not strict but are isomorphic to strict objects, one makes obvious adjustments.
\end{proof}

The main tool in proving controlled
excision in the boundedly controlled
$K$-theory will be the following localization theorem.

\begin{ThmRefName}{Schlichting}{Theorem 2.1 of Schlichting \cite{mS:03}}
Let $\Z$ be an idempotent complete right s-filtering subcategory
of an exact category $\B$.  Then the sequence of exact categories
$\Z \to \B \to {\B}/{\Z}$ induces a homotopy fibration of K-theory
spectra
\[
K(\Z) \longrightarrow K(\B) \longrightarrow K({\B}/{\Z}).
\]
\end{ThmRefName}

\SecRef{Bounded excision theorem}{LocK}

The proof of controlled excision in the boundedly controlled
$G$-theory requires the context of
Waldhausen $K$-theory of derived
categories.

\begin{DefRefName}{Wald}{Waldhausen categories}
A \define{Waldhausen category} is a category $\D$ with a zero
object $0$ together with two chosen subcategories of
\define{cofibrations} $\co (\D)$ and \define{weak equivalences}
$\w (\D)$ satisfying the four axioms:
\begin{enumerate}
\item every isomorphism in $\D$ is in both $\co (\D)$ and $\w (\D)$,
\item every map $0 \to D$ in $\D$ is in $\co (\D)$,
\item if $A \to B \in \co (\D)$ and $A \to C \in \D$ then
the pushout $B \cup_A C$ exists in $\D$, and the canonical map $C
\to B \cup_A C$ is in $\co (\D)$,
\item (``gluing lemma'') given a commutative diagram
\[
\begin{CD}
B @<{a}<< A @>>> C \\
@VVV @VVV @VVV \\
B' @<{a'}<< A' @>>> C'
\end{CD}
\]
in $\D$, where the morphisms $a$ and $a'$ are in $\co (\D)$ and
the vertical maps are in $\w (\D)$, the induced map $B \cup_A C
\to B' \cup_{A'} C'$ is also in $\w (\D)$.
\end{enumerate}
A Waldhausen category $\D$ with weak equivalences $\w (\D)$ is
often denoted by $\wD$ as a reminder of the choice. A functor
between Waldhausen categories is exact if it preserves
cofibrations and weak equivalences.
\end{DefRefName}

A Waldhausen category may or may not satisfy the following
additional axioms.
\begin{SatAx}
Given two morphisms $\phi \colon F \to G$ and $\psi \colon G \to
H$ in $\D$, if any two of $\phi$, $\psi$, or $\psi  \phi$, are in
$\w (\D)$ then so is the third.
\end{SatAx}

\begin{ExtAx}
Given a commutative diagram
\[
\begin{CD}
 F @>>> G @>>> H\\
 @VV{\phi}V @VV{\psi}V @VV{\mu}V\\
 F' @>>> G' @>>> H'
\end{CD}
\]
with exact rows, if both $\phi$ and $\mu$ are in $\w (\D)$ then so
is $\psi$.
\end{ExtAx}

A \textit{cylinder functor} on $\D$ is a functor $C$ from the
category of morphisms $f \colon F \to G$ in $\D$ to $\D$ together
with three natural transformations $j_1 \colon F \to C(f)$, $j_2
\colon G \to C(f)$, and $p \colon C(f) \to G$ which satisfy $p j_2 =
\id_G$ and $p j_1 = f$ for all $f$.
The cylinder functor has to satisfy a number of
properties listed in point 1.3.1 of \cite{rTtT:90}.
They will be
rather automatic for the functors we construct later.

\begin{CylAx}
A cylinder functor $C$ satisfies this axiom if for all morphisms
$f \colon F \to G$ the required map $p$ is in $\w (\D)$.
\end{CylAx}

Let $\D$ be a small Waldhausen category with respect to two
categories of weak equivalences $\bsv(\D) \subset \w(\D)$ with a
cylinder functor $T$ both for $\vD$ and for $\wD$ satisfying the
cylinder axiom for $\wD$. Suppose also that $\w(\D)$ satisfies the
extension and saturation axioms.

Define $\vDw$ to be the full
subcategory of $\vD$ whose objects are $F$ such that $0 \to F \in
\w(\D)$. Then $\vDw$ is a small Waldhausen category with
cofibrations $\co (\Dw) = \co (\D) \cap \Dw$ and weak
equivalences $\bsv (\Dw) = \bsv (\D) \cap \Dw$. The cylinder
functor $T$ for $\vD$ induces a cylinder functor for $\vDw$.  If
$T$ satisfies the cylinder axiom then the induced functor does so
too.

\begin{ThmRefName}{ApprThm}{Approximation theorem}
Let $E \colon \D_1 \rightarrow \D_2$ be an exact functor between
two small saturated Waldhausen categories. It induces a map of
$K$-theory spectra
\[
  K (E) \colon K (\D_1) \longrightarrow K (\D_2).
\]
Assume that $\D_1$ has a cylinder functor satisfying the cylinder
axiom. If $E$ satisfies two conditions:
\begin{enumerate}
\item a morphism $f \in \D_1$ is in $\w(\D_1)$
if and only if $E (f) \in \D_2$ is in $\w(\D_2)$,
\item for any object $D_1 \in \D_1$ and any morphism
$g \colon E(D_1) \to D_2$ in $\D_2$, there is an object $D'_1 \in
\D_1$, a morphism $f \colon D_1 \to D'_1$ in $\D_1$, and a weak
equivalence $g' \colon E(D'_1) \rightarrow D_2 \in \w(\D_2)$ such
that $g = g' E(f)$,
\end{enumerate}
then $K (E)$ is a homotopy equivalence.
\end{ThmRefName}

\begin{proof}
This is Theorem~1.6.7 of \cite{fW:85}. The presence of the
cylinder functor with the cylinder axiom allows to make condition
2 weaker than that of Waldhausen, see point 1.9.1 in
\cite{rTtT:90}.
\end{proof}

\begin{DefRef}{PFff}
In any additive category, a sequence of morphisms
\[
E^{\subdot} \colon \quad 0 \longrightarrow E^1 \xrightarrow{\ d_1
\ } E^2 \xrightarrow{\ d_2 \ }\ \dots\ \xrightarrow{\ d_{n-1} \ }
E^n \longrightarrow 0
\]
is called a \textit{(bounded) chain complex} if the compositions
$d_{i+1} d_i$ are the zero maps for all $i = 1$,\dots, $n-1$.  A
\textit{chain map} $f \colon F^{\subdot} \to E^{\subdot}$ is a
collection of morphisms $f^i \colon F^i \to E^i$ such that $f^i
d_i = d_i f^i$.

A chain map $f$ is \textit{null-homotopic} if
there are morphisms $s_i \colon F^{i+1} \to E^i$ such that $f = ds
+ sd$.  Two chain maps $f$, $g \colon F^{\subdot} \to E^{\subdot}$
are \textit{chain homotopic} if $f-g$ is null-homotopic. Now $f$
is a \textit{chain homotopy equivalence} if there is a chain map
$h \colon E^i \to F^i$ such that the compositions $fh$ and $hf$
are chain homotopic to the respective
identity maps.
\end{DefRef}

The Waldhausen structures on categories of bounded chain complexes
are based on homotopy equivalence as a weakening of the notion of
isomorphism of chain complexes.

\begin{DefRef}{PF}
A sequence of maps in an exact category is called \textit{acyclic}
if it is assembled out of short exact sequences in the sense that
each map factors as the composition of the cokernel of the
preceding map and the kernel of the succeeding map.
\end{DefRef}

It is known that the class of acyclic complexes in an exact
category is closed under isomorphisms in the homotopy category if
and only if the category is idempotent complete, which is also
equivalent to the property that each contractible chain complex is
acyclic, cf.\ \cite[sec.\ 11]{bK:96}.

\begin{DefRef}{PrepFil23}
Given an exact category $\E$, there is a standard choice for the
Waldhausen structure on the derived category $\E'$ of bounded
chain complexes in $\E$.
The cofibrations $\mathrm{\textbf{co}} (\E')$ are the degree-wise admissible
monomorphisms.
The
weak equivalences $\boldsymbol{v}(\E')$
are the chain maps whose
mapping cones are homotopy equivalent to acyclic complexes.

We will denote this Waldhausen structure by $\vE'$.
\end{DefRef}

\begin{PropRef}{FibThApp2}
The category $\vE'$ is a Waldhausen category satisfying the
extension and saturation axioms and has cylinder functor
satisfying the cylinder axiom.
\end{PropRef}

\begin{proof}
The pushouts along cofibrations in $\E'$ are the complexes of
pushouts in each degree.  All standard Waldhausen axioms including
the gluing lemma are clearly satisfied. The saturation and the
extension axioms are also clear.

The cylinder functor $C$ for
$\vE'$ is defined using the canonical homotopy pushout as in point
1.1.2 in Thomason--Trobaugh \cite{rTtT:90}. Given a chain map $f
\colon F \to G$, $C(f)$ is the canonical homotopy pushout of $f$
and the identity $\id \colon F \to F$. With this construction, the
map $p \colon C(f) \to G$ is a chain homotopy equivalence, so the
cylinder axiom is also satisfied.
\end{proof}

\begin{DefRef}{PrepFil2}
There are two choices for the Waldhausen structure on the bounded
derived category $\B'$ for the exact boundedly controlled category
$\B = \B (X, \E)$. One is $\vB'$ as in
Definition \refD{PrepFil23}. Given a metric subspace $Z$ in $X$,
the other choice for the weak equivalences $\w (\B')$ is the chain
maps whose mapping cones are homotopy equivalent to acyclic
complexes in the quotient $\mathbf{B}/\mathbf{Z}$.

We will denote the second Waldhausen structure by $\wB'$.
\end{DefRef}

\begin{CorRef}{FibSt}
The categories $\vB'$ and $\wB'$ are Waldhausen categories
satisfying the extension and saturation axioms and have cylinder
functors satisfying the cylinder axiom.
\end{CorRef}

\begin{proof}
All axioms and constructions, including the cylinder functor, for
$\wB'$ are inherited from $\vB'$.
\end{proof}

The $K$-theory functor from the category of small Waldhausen
categories $\D$ and exact functors to connective spectra is
defined in terms of $S_{\subdot}$-construction as in Waldhausen
\cite{fW:85}. It extends to simplicial categories $\D$ with
cofibrations and weak equivalences and inductively gives the
connective spectrum
\[
n \mapsto \vert \bfw S_{\subdot}^{(n)} \D
\vert.
\]
We obtain the functor assigning to $\D$ the connective
$\Omega$-spectrum
\[
K (\D) = \Omega^{\infty} \vert \bfw S_{\subdot}^{(\infty)} \D
\vert = \colim{n \ge 1} \Omega^n \vert \bfw S_{\subdot}^{(n)} \D
\vert
\]
representing the Waldhausen algebraic $K$-theory of $\D$. For
example, if $\D$ is the additive category of free finitely
generated $R$-modules with the canonical Waldhausen structure,
then the stable homotopy groups of $K (\D)$ are the usual
$K$-groups of the ring $R$.  In fact, there is a general
identification of the two theories.

Recall that for any exact
category $\E$, the derived category $\E'$ has the Waldhausen
structure $\vE'$ as in Definition \refD{PrepFil23}.

\begin{ThmRef}{Same}
The Quillen $K$-theory of an exact category $\E$ is equivalent to
the Waldhausen $K$-theory of $\vE'$.
\end{ThmRef}

\begin{proof}
The proof is based on repeated applications of the additivity
theorem, cf. Thomason's Theorem 1.11.7 \cite{rTtT:90}. Thomason's
proof of his Theorem 1.11.7 can be repeated verbatim here.  It is
in fact simpler in this case since condition 1.11.3.1 is not
required.
\end{proof}

Let $\E$ be a Grothendieck subcategory of a cocomplete abelian
category $\A$ and let $(Z,X)$ be a pair of proper metric spaces.

We will use the notation $\B = \B (X, \E)$ and $\Z = \B (X,
\E)_{<Z}$.

\begin{ThmRefName}{LocText}{Localization}
If $\E$ is idempotent complete, there is a homotopy fibration
\[
K (Z, \E) \longrightarrow K (X,\E) \longrightarrow K ({\B}/{\Z}).
\]
\end{ThmRefName}

This is a direct consequence of Theorem \refT{Schlichting}
as soon as we identify $K (Z, \E)$ with $K (\Z) = K (X,
\E)_{<Z}$.

Recall that the \textit{essential full image} of a functor $F
\colon \C \to \D$ is the full subcategory of $\D$ whose objects
are those $D$ such that $D \cong F(C)$ for some $C$ from $\C$.

There is a fully faithful embedding
\[
\epsilon \colon \B (Z,\E) \longrightarrow
\B (X,\E)
\]
defined by associating to each filtered object $F \in \B
(Z,\E)$ the extension $\epsilon (F) \in \B (X,\E)$ given by
\[
\epsilon (F)(S) = F (S \cap Z).
\]
It is clear that $\epsilon (F)$ is strict for strict $F$.

\begin{LemRef}{LemNeed}
The essential full image of $\B (Z,\E)$ in $\B (X,\E)$ is the
Grothendieck subcategory $\B (X,\E)_{<Z}$.
\end{LemRef}

\begin{proof}
Of course for each $F$ in $\B (Z,\E)$, the image $\epsilon (F)$ is in $\B (X,\E)_{<Z}$.  Now if $G(X) =
G(Z[d])$ is an object of $\B (X,\E)_{<Z}$ then there is a bounded function $r
\colon Z[d] \to Z$, bounded by $d$, which gives an object $R = R(G)$ of $\B (Z,\E)$
by the assignment $R(S) = G(r^{-1}(S))$.

If $G$ is strict then the new object $R$ is
$\E$-local and strict with $\ell_R = \ell_G + d$. Since the identity map
$\id \colon R \to G$ is boundedly bicontrolled with $\fil (\id) \le 2d$, it is an
isomorphism in $\B (X,\E)$.
\end{proof}

\begin{CorRef}{CompReal}
For any pair of proper metric spaces $Z \subset X$, there is a
weak equivalence $K (Z,\E) \simeq K (X,\E)_{<Z}$.
\end{CorRef}

Now Theorem \refT{LocText} follows from the localization fibration in Theorem \refT{Schlichting}.

\bigskip

The results of Theorem \refT{LocText} and Lemma \refL{LemNeed} can be generalized to a more general and convenient geometric situation.

Suppose $Z$ is an arbitrary subset of a proper metric space $X$.  It is a metric subspace with the metric which is the restriction of the metric in $X$.  When $Z$ is a closed subset then the closed metric balls in $Z$ are closed subsets of closed metric balls in $X$ and, therefore, compact.
If $Z$  is an arbitrary subset of $X$, the closure $\overline{Z}$ is a proper metric subspace.
There is an inclusion $\overline{Z} \subset Z[\epsilon]$ for any $\epsilon >0$, so there is
an $\epsilon$-bounded retraction of $\overline{Z}$ onto $Z$.
In addition to the obvious equivalence of categories $\B (X, \E)_{<\overline{Z}} = \B (X, \E)_{<Z}$, the retraction also induces an isomorphism $\B (\overline{Z}, \E) \cong \B (Z, \E)$.

\begin{CorRefName}{LocArb}{Localization}
Suppose $X$ is a proper metric space and $Z$ is a subset with the induced metric.
There is a
weak homotopy equivalence
\[
K (\overline{Z}, \E) \simeq K (Z, \E).
\]
If $\E$ is idempotent complete, there is a homotopy fibration
\[
K (Z, \E) \longrightarrow K (X,\E) \longrightarrow K ({\B}/{\Z}).
\]
\end{CorRefName}

\bigskip

The computational tools from
nonconnective bounded $K$-theory, the controlled
excision
theorems \cite{mCeP:97, ePcW:85,
ePcW:89}, can now be adapted to
$\B(X,\E)$. We will obtain a
direct analogue, which is one of the main
results of this paper.

Let $\E$ be a Grothendieck subcategory of a cocomplete abelian
category $\A$ and let $X$ be a proper metric space. Suppose $X_1$
and $X_2$ are subspaces in a proper metric space $X$, and $X = X_1
\cup X_2$.

We use the notation $\B = \B (X,\E)$, $\B_i= \B
(X,\E)_{<X_i}$ for $i=1$ or $2$, and $\B_{12}$ for the intersection
$\B_1 \cap \B_2$.

Now there is a commutative diagram
\[
\begin{CD}
K (\B_{12}) @>>> K (\B_1) @>>> K ({\B_1}/{\B_{12}}) \\
@VVV @VVV @VV{K(I)}V \\
K (\B_2) @>>> K (\B) @>>> K ({\B}/{\B_2})
\end{CD} \tag{$\dagger$}
\]
where the rows are homotopy fibrations from Theorem \refT{Schlichting} and $I \colon {\B_1}/{\B_{12}} \to {\B}/{\B_2}$ is the exact functor induced from the exact inclusion $I \colon \B_1 \to \B$.
We observe that $I$ is not necessarily full and,
therefore, not an isomorphism of categories as in similar
applications in \cite{mCeP:97} and \cite{rS:89}.
Nevertheless, we claim that
$K (I)$ is a weak equivalence.

\begin{LemRef}{CharWE}
Let $Z$ be a subset of $X$, so $\Z = \B (X, \E)_{<Z}$ is a
Grothendieck subcategory of $\B$.
If $f^{\subdot}$ is  a
degreewise admissible monomorphism with cokernels in $\Z$ then
$f^{\subdot}$ is a weak equivalence in $\boldsymbol{v}
({\B}/{\Z})'$.
\end{LemRef}

\begin{proof}
The mapping cone $Cf^{\subdot}$ of $f^{\subdot}$ is homotopy equivalent to the cokernel of
$f^{\subdot}$, which is an acyclic complex in
${\B}/{\Z}$.
\end{proof}

\begin{LemRef}{IdTarget}
The map
\[
K(I) \colon K  ({\B_1}/{\B_{12}}) \longrightarrow K ({\B}/{\B_2})
\]
is a weak equivalence.
\end{LemRef}

\begin{proof}
We will apply the Approximation Theorem to $I$.
The first condition is clear.
To check the second condition, consider
\[
F^{\subdot} \colon \quad 0 \longrightarrow F^1 \xrightarrow{\
\phi_1\ } F^2 \xrightarrow{\ \phi_2\ }\ \dots\
\xrightarrow{\ \phi_{n-1}\ } F^n \longrightarrow 0
\]
in $\B_1$ and a chain map $g \colon F^{\subdot} \to G^{\subdot}$ to some complex
\[
G^{\subdot} \colon \quad 0 \longrightarrow G^1 \xrightarrow{\
\psi_1\ } G^2 \xrightarrow{\ \psi_2\ }\ \dots\
\xrightarrow{\ \psi_{n-1}\ } G^n \longrightarrow 0
\]
in $\B$.
Without loss of generality, suppose that each $G^i$ is $D$-lean, suppose that $F^i = F^i (X_1[K])$ for all $i$, and suppose that $b$ serves as a bound for all $\phi_i$, $\psi_i$, and $g_i$.
We define
\[
F^{\prime{i}} = G^{{i}} (X_1[K + D + 3ib])
\]
and define $\xi_{i} \colon F^{\prime{i}} \to F^{\prime{i+1}}$ to be the restrictions of
$\psi_{i}$ to $F^{\prime{i}}$. This gives a chain subcomplex
$(F^{\prime{i}},\xi_{i})$ of $(G^{i},\psi_{i})$ in $\B_1$ with the inclusion $s \colon F^{\prime{i}} \to G^{i}$.
Notice that the choices give the induced chain map $\overline{g} \colon F^{\subdot} \to F^{\prime\subdot}$ in $\B_1$ so that $g = s \circ I(\overline{g})$.

We will argue that $C^{\subdot} = \coker (s)$ is in $\B_2$.  Given that, $s$ is a weak equivalence in
$\boldsymbol{v}
({\B}/{\B_2})'$ by Lemma \refL{CharWE}.

For each $x \in X_1$ and each $i$, we have $G^i (B_D (x)) \subset F^{\prime{i}}$, so $C^i (B_D (x)) = 0$.  By Proposition \refP{lninpres}, part (3b), $C^i$ is $D$-lean, therefore
\[
C^i (X) = \sum_{x \in X} C^i (B_D (x)) = \sum_{x \in X \backslash X_1} C^i (B_D (x)) \subset C^i ((X \backslash X_1)[D]) \subset C^i (X_2).
\]
So the complex $C^{\subdot}$ is indeed in $\B_2$.
\end{proof}

Let $\mathbb{Z}$, $\mathbb{Z}^{\ge 0}$, and $\mathbb{Z}^{\le 0}$
denote the metric spaces of integers, nonnegative integers, and
nonpositive integers with the restriction of the usual metric on the
real line $\mathbb{R}$.  Let $\E$ be an idempotent complete
Grothendieck category of an abelian category $\A$. Then for any
proper metric space $X$, we have the following instance of
commutative diagram ($\dagger$)

\[
\begin{CD}
K (X, \E) @>>> K (X \times \mathbb{Z}^{\ge 0}, \E) @>>> K ({\B_1}/{\B_{12}}) \\
@VVV @VVV @VV{K(I)}V \\
K (X \times \mathbb{Z}^{\le 0}, \E) @>>> K (X \times \mathbb{Z}, \E) @>>> K ({\B}/{\B_2})
\end{CD}
\]

\begin{LemRef}{ESw}
The spectra $K (X \times \mathbb{Z}^{\ge 0}, \E)$ and $K (X \times \mathbb{Z}^{\le 0}, \E)$ are contractible.
\end{LemRef}

\begin{proof}
This follows from the fact that these controlled categories are
flasque, that is, the usual shift functor $T$ in the positive
(respectively negative)
direction along $\mathbb{Z}^{\ge 0}$ (respectively $\mathbb{Z}^{\le 0}$) interpreted in the obvious
way is an exact endofunctor, and there is a natural equivalence $1
\oplus \pm T \cong \pm T$. Contractibility follows from the
additivity theorem, cf.\ Pedersen--Weibel \cite{ePcW:85}.
\end{proof}

In view of Lemma \refL{IdTarget}, we
obtain a map
\[
K (X, \E) \longrightarrow \Omega K (X \times \mathbb{Z}, \mathbf{E})
\]
which induces isomorphisms of $K$-groups
in positive dimensions.
Iterations of this
construction give weak equivalences
\[
\Omega^k K (X \times \mathbb{Z}^k,
\mathbf{E})
\longrightarrow \Omega ^{k+1}
K (X \times \mathbb{Z}^{k+1}, \E)
\]
for $k \ge 2$.

\begin{DefRefName}{RealExQ}{Nonconnective controlled $K$-theory}
The \textit{nonconnective controlled $K$-theory} of $\E$, relative
to the embedding $\epsilon \colon \E \to \A$, over a proper metric
space $X$ is the spectrum
\[
\Knc_{\epsilon} (X, \E) \overset{ \text{def} }{=} \hocolim{k}
\Omega^{k} K (X \times \mathbb{Z}^k, \E).
\]
\end{DefRefName}

Since $\B (X, \E)$ can be identified with $\E$ for a bounded metric
space $X$, this definition gives the \textit{nonconnective
$K$-theory} of $\E$
\[
\Knc_{\epsilon} (\E) \overset{ \text{def} }{=} \hocolim{k > 0}
\Omega^{k} K (\mathbb{Z}^k, \E).
\]

As $\Knc_{\epsilon} (\E)$ is an $\Omega$-spectrum in positive
dimensions, the positive homotopy groups of $\Knc (\E)$ coincide
with those of $K (\E)$, as desired.  The class group
$K_{\epsilon,0} (\E)$ is the class group of the idempotent
completion $K_0 (\E\sphat \,)$.

The first known delooping of the $K$-theory of a general exact
category with these properties is due to M.~Schlichting
\cite{mS:99}, however the construction here is different and is
required in the excision theorem ahead.

\begin{ExRef}{CaDe}
If $\E$ is an arbitrary small exact category, there is the full
Gabriel--Quillen embedding of $\E$ in the cocomplete abelian
category $\A$ of left exact functors $\E^{\textrm{op}} \to \Mod
(\mathbb{Z})$ with the standard exact structure. The embedding is
always closed under extensions in $\A$.  It is not necessarily a
Grothendieck subcategory, as when $\E$ is not balanced. But
if $\E$ is abelian, for example, this gives a canonical delooping of $K (\E)$.
\end{ExRef}

\begin{ExRef}{Gtheory}
One may start with the cocomplete abelian category $\Mod (R)$ of
modules over a ring $R$ with the standard abelian exact structure
where the admissible monomorphisms and epimorphisms are
respectively all monics and epis.  If $R$ is a Noetherian ring,
the subcategory $\E$ may be taken to be the noncocomplete abelian
category of finitely generated $R$-modules $\Modf(R)$.  Now $\Knc
(\E)$ gives the algebraic $G$-theory of $R$ which we denote by
$\Gnc (R)$.  We also use notation $\B (X, R)$ for $\B (X, \E)$.
\end{ExRef}

\begin{DefRefName}{Gth}{Nonconnective controlled $G$-theory}
The \textit{nonconnective controlled $G$-theory} of
$X$-filtered modules over $R$ is defined as
\[
\Gnc (X,R) = \Knc \B (X, \E).
\]
\end{DefRefName}

\begin{ExRef}{CocE}
The negative $K$-theory of a regular ring $R$ is trivial in the
sense that $\K_i (\Modf(R)) = 0$ for all $i < 0$. This is
well-known in Bass' theory \cite{hB:68}.  A proof that the
negative $K$-theory is trivial for general abelian categories can
be given using the same strategy as in chapter 9 of \cite{mS:99}.
\end{ExRef}

Of course, when the exact category $\E$ is itself cocomplete, its
$K$-theory is contractible because of the Eilenberg swindle type
argument.

We finally prove the main result of this section.

\begin{ThmRefName}{Exc}{Nonconnective excision}
Let $\E$ be a Grothendieck subcategory of a cocomplete abelian
category $\A$ and let $X$ be a proper metric space. Suppose $X_1$
and $X_2$ are subsets of $X$, and $X = X_1
\cup X_2$.
Using the notation $\B = \B (X,\E)$, $\B_i= \B
(X,\E)_{<X_i}$ for $i=1$ or $2$, and $\B_{12}$ for the intersection
$\B_1 \cap \B_2$,
there is a homotopy pushout diagram of spectra
\[
\begin{CD}
\Knc (\B_{12}) @>>> \Knc (\B_{1}) \\
@VVV @VVV \\
\Knc (\B_{2}) @>>> \Knc (\B)
\end{CD}
\]
where the maps of spectra are induced from the exact inclusions.
\end{ThmRefName}

\begin{proof}
Let us write $S^k \B$ for $\B (X \times \mathbb{Z}^k, \E)$
whenever $\B$ is the boundedly controlled category for a general
metric space $X$.  If $\Z$ is a subset of $X$, consider the
fibration
\[
K (Z, \E) \longrightarrow K (X,\E) \longrightarrow K ({\B}/{\Z})
\]
from Corollary \refC{LocArb}.  Notice that there is a map
\[
K ({\B}/{\Z}) \longrightarrow \Omega K ({S\B}/{S\Z})
\]
which is a weak equivalence in positive dimensions by the Five
Lemma. If one defines
\[
\Knc ({\B}/{\Z}) = \hocolim{k}
\Omega^{k} K ({S^k \B}/{S^k \Z}),
\]
there is an induced fibration
\[
\Knc (Z, \E) \longrightarrow \Knc (X,\E) \longrightarrow \Knc ({\B}/{\Z})
\]
The theorem follows from the commutative diagram
\[
\begin{CD}
\Knc (\B_{12}) @>>> \Knc (\B_1) @>>> \Knc ({\B_1}/{\B_{12}}) \\
@VVV @VVV @VV{\Knc (I)}V \\
\Knc (\B_2) @>>> \Knc (\B) @>>> \Knc ({\B}/{\B_2})
\end{CD}
\]
and the fact that now
\[
\Knc (I) \colon \Knc  ({\B_1}/{\B_{12}})
\longrightarrow \Knc ({\B}/{\B_2})
\]
is a weak equivalence.
\end{proof}

\begin{RemRef}{NoConnExc}
As in other versions of controlled $K$-theory, there is no
excision theorem similar to Theorem \refT{Exc} which employs the
connective $K$-theory.
\end{RemRef}

\SecRef{Equivariant theory and the Novikov conjecture}{FuncOP}

First we establish functoriality properties of the bounded $G$-theory.

\begin{DefRef}{POIUY}
A map $f \colon X \to Y$
of metric spaces $f$ is \define{quasi-bi-Lipschitz} if there is a real positive function $l$ such that
\begin{gather}
\dist(x_1, x_2) \le r \ \Longrightarrow
\dist (f(x_1), f(x_2)) \le l(r), \notag \\
\dist (f(x_1), f(x_2)) \le r \ \Longrightarrow \ \dist(x_1, x_2)
\le l(r) \notag
\end{gather}
for all $x_1$, $x_2 \in X$.

Any bounded function $f \colon X \to X$, with the property that
$\dist(x, f(x)) \le D$ for all $x \in X$ and a fixed number $D \ge 0$, is
quasi-bi-Lipschitz with $l(r)=r+2D$. An isometry $g \colon X \to Y$
is quasi-bi-Lipschitz with $l (r) = r$.
If only the first of the two
conditions is satisfied, the map $f$ is called
\define{bornological}.

We will say that $f \colon X \to Y$
is a \define{quasi-bi-Lipschitz equivalence} if there is a map $g \colon Y \to X$ so that
both $f$ and $g$ are quasi-bi-Lipschitz and the compositions $f \circ g$ and $g \circ f$ are bounded maps.
\end{DefRef}

\begin{PropRef}{FunctB}
Consider the category of proper metric spaces $X$ and
quasi-bi-Lipschitz maps and the category of Noetherian rings
$R$.
Then $\B (X,R)$ is a bifunctor covariant in the first
variable and contravariant in the second variable to exact
categories and exact functors. Composing with the covariant
functor $\Knc$ from Example \refE{Gtheory} gives the spectrum-valued bifunctor
$\Gnc (X,R)$.
\end{PropRef}

\begin{proof}
If $f \colon X \to Y$ is a quasi-bi-Lipschitz map, the
functor
\[
f_{\ast} \colon \B (X,R) \longrightarrow \B(Y,R)
\]
is given on objects by
\[
f_{\ast} F(S) = F(f^{-1}(S)).
\]
Using the containment
\[
f^{-1}(S)[D] \subset f^{-1}(S[l(D)]),
\]
one sees that if $\phi \in
\B (X,R)$ is a boundedly bicontrolled morphism with $\fil (\phi)
\le D$ then $f_{\ast} \phi$ is boundedly bicontrolled with $\fil
(f_{\ast} \phi) \le l(D)$.
\end{proof}

A subset $W$ of a metric space $X$ is \define{boundedly dense} or
\define{commensurable} if $W[D] = X$ for some $D \ge 0$.

\begin{PropRef}{IsomOK}
For a commensurable metric subspace $W$ of $X$, there is a natural
exact equivalence of categories
\[
i \colon \B (W,R) \longrightarrow \B (X,R)
\]
and the
induced weak homotopy equivalence
\[
i_{\ast} \colon \Gnc (W,R) \longrightarrow \Gnc (X,R).
\]
\end{PropRef}

\begin{proof}
Any surjective quasi-bi-Lipschitz equivalence $f \colon X \to Y$
induces two functors on filtered modules. One is contravariant
\[
f^{\ast} \colon \B (Y,R) \longrightarrow \B (X,R)
\]
given by $f^{\ast}F(S) =
F(f(S))$; the other is covariant
\[
f_{\ast} \colon \B (X,R) \longrightarrow \B
(Y,R)
\]
given by $f_{\ast}F(S) = F(f^{-1}(S))$, so that $f^{\ast}
f_{\ast} = \id$.

Even when $f$ is not surjective, there is the
endofunctor $\omega = f^{-1} f$ of $\mathcal{P}(X)$ which induces
an endofunctor $\omega_{\ast}$ of $\B (X,R)$. If $f \colon X \to
X$ is bounded, that is $d(x,f(x)) \le D$ for some $D \ge 0$ and
all $x \in X$, there is always an isomorphism $\omega_{\ast}(F)
\cong F$ induced by the identity on $F(X)$. This shows that
$f_{\ast} F \cong F$ for all $F \in \B (X,R)$.

If $W \subset
X$ is commensurable, there is a bounded surjection $f \colon X \to
W$, so $f$ induces a natural transformation $\eta \colon \id \to
f_{\ast}$ where all $\eta(F)$ are isomorphisms.
\end{proof}

\begin{CorRef}{BddMod}
If $X$ is a bounded metric space then the natural equivalence
\[
\B (X,R) \cong \B (\point,R) = \Modf(R)
\]
induces a weak equivalence $\Gnc (X,R) \simeq \Gnc (R)$ on
the level of $K$-theory.
\end{CorRef}

Given a group $\Gamma$ with a left action on $X$ by
quasi-bi-Lipschitz equivalences, there is a natural action of
$\Gamma$ on $\B (X,R)$ induced from the action on the power set
$\mathcal{P}(X)$. However, this is not the correct choice for a
useful equivariant controlled theory for essentially the same
reasons as in the discussion of geometric modules in
\cite[ch.\,VI]{gC:95}.

\begin{DefRef}{EqCatDef}
Let $\EGamma$ be the category with the object set $\Gamma$ and the
unique morphism $\mu \colon \gamma_1 \to \gamma_2$ for any pair
$\gamma_1$, $\gamma_2 \in \Gamma$. There is a left $\Gamma$-action
on $\EGamma$ induced by the left multiplication in $\Gamma$.

If $\mathcal{C}$ is a small category with left $\Gamma$-action,
then the category of functors
$\mathcal{C}_{\Gamma}=\Fun(\EGamma,\mathcal{C})$ is another
category with the $\Gamma$-action given on objects by
\[
\gamma(F)(\gamma')=\gamma F (\gamma^{-1} \gamma')
\]
and
\[
\gamma(F)(\mu)=\gamma F (\gamma^{-1} \mu).
\]
It is always
nonequivariantly equivalent to $\mathcal{C}$. The fixed
subcategory $\Fun(\EGamma,\mathcal{C})^{\Gamma} \subset
\mathcal{C}_{\Gamma}$ consists of equivariant functors and
equivariant natural transformations.

Explicitly, when $C = \B (X,R)$ with the $\Gamma$-action described
above, the objects of $\B_{\Gamma} (X,R)^{\Gamma}$ are the pairs
$(F,\psi)$ where $F \in \B (X,R)$ and $\psi$ is a function on
$\Gamma$ with $\psi (\gamma) \in \Hom (F,\gamma F)$ such that
\[
\psi(1) = 1 \quad \mathrm{and} \quad \psi (\gamma_1 \gamma_2) =
\gamma_1 \psi(\gamma_2)  \psi (\gamma_1).
\]
These conditions imply that $\psi (\gamma)$ is always an
isomorphism as in \cite{rT:83}. The set of morphisms $(F,\psi) \to
(F',\psi')$ consists of the morphisms $\phi \colon F \to F'$ in
$\B (X,R)$ such that the squares
\[
\begin{CD}
F @>{\psi (\gamma)}>> \gamma F \\
@V{\phi}VV @VV{\gamma \phi}V \\
F' @>{\psi' (\gamma)}>> \gamma F'
\end{CD}
\]
commute for all $\gamma \in \Gamma$. A slightly more refined
theory is obtained by replacing $\B_{\Gamma} (X,R)$ with the full
subcategory $\B_{\Gamma, 0} (X,R)$ of functors sending all
morphisms of $\EGamma$ to filtration $0$ maps. So $\B_{\Gamma, 0}
(X,R)^{\Gamma}$ consists of $(F, \psi)$ with $\fil \psi (\gamma) =
0$ for all $\gamma \in \Gamma$.
\end{DefRef}

\begin{PropRef}{ExStr}
The fixed point category $\B_{\Gamma, 0} (X,R)^{\Gamma}$ is exact.
\end{PropRef}

\begin{proof}
The exact structure is inherited from $\B (X,R)$ in the sense that
a morphism $\phi \colon (F, \psi) \to (F', \psi')$ is an
admissible monomorphism or epimorphism if the map $\phi \colon F
\to F'$ is in $\mB (X,R)$ or $\eB (X,R)$ respectively. The fact
that this is an exact structure follows from the proof of Theorem
\refT{UbisWAb} by observing that all constructions in that proof
produce equivariant objects and morphisms.
\end{proof}

\begin{RemRef}{UIYOU}
One exact structure in the category of finitely generated
$R[\Gamma]$-modules $\Modf (R[\Gamma])$ for a Noetherian ring $R$
consists of short exact sequences of $R[\Gamma]$-modules with
finitely generated kernels and quotients.
When $R[\Gamma]$ is Noetherian, so that $\Modf
(R[\Gamma])$ is an abelian category, this coincides with the
conventional choice of all injections for admissible monomorphisms
and all surjections for admissible epimorphisms.

However, there is
a reasonable conjecture of P.~Hall that only polycyclic-by-finite
groups have Noetherian group rings, cf.\ Question 32 in Farkas
\cite{dF:80}.
\end{RemRef}

We are going to define a new exact structure on a subcategory $\B
(R[\Gamma])$ of $\Modf (R[\Gamma])$ and relate it to the exact
category $\B_{\Gamma,0} (X,R)^{\Gamma}$.

\begin{DefRef}{QWCFRE}
The \textit{word metric} on a finitely generated group
$\Gamma$ with a fixed generating set $\Omega$ is the path metric
induced from the condition that $\dist (\gamma, \omega \gamma) =
1$ whenever $\gamma \in \Gamma$ and $\omega \in \Omega$.
\end{DefRef}

This metric
clearly makes $\Gamma$ a proper metric space.  We will use the
notation $B_d (\gamma)$ for the metric ball of radius $d$ centered
at $\gamma$.

\begin{DefRef}{AndBack}
Given a finitely generated $R[\Gamma]$-module $F$, fix a finite
generating set $\Sigma$ for $F$ and define a $\Gamma$-filtration of the
$R$-module $F$ by
\[
F(S) = \langle S \Sigma \rangle_R,
\]
the
$R$-submodule of $F$ generated by $S \Sigma$.  Let $s(F, \Sigma)$
stand for the resulting $\Gamma$-filtered $R$-module.
\end{DefRef}

\begin{LemRef}{AlmPull}
Every $R[\Gamma]$-homomorphism
\[
\phi \colon F \longrightarrow G
\]
between finitely generated modules
is boundedly controlled as an
$R$-homomorphism between the filtered $R$-modules
\[
\phi \colon s (F, \Sigma_F)
\longrightarrow s (G, \Sigma_G)
\]
with respect to any choice of the finite
generating sets $\Sigma_F$ and $\Sigma_G$.
\end{LemRef}

\begin{proof}
Consider $x \in F(S) = \langle S \Sigma_F \rangle_R$, then
\[
x = \sum_{s, \sigma} r_{s, \sigma} s \sigma
\]
for a finite collection of pairs $s \in S$, $\sigma \in \Sigma_F$.
Since $F (\{ e \}) = \langle \Sigma_F \rangle_R$ for the identity
element $e$ in $\Gamma$, there is a number $d \ge 0$ such that
\[
\phi F (\{ e \}) \subset G(B_d (e)).
\]
Therefore,
\[
\phi (x) = \sum_{s, \sigma} r_{s, \sigma} \phi (s \sigma) =
\sum_{s, \sigma} r_{s, \sigma} s \phi (\sigma) \subset \sum_{s \in
S} s G(B_d (e)) \subset G (S[d])
\]
because the left translation action by any element $s \in S$ on $B_d (e)$
in $\Gamma$ is an isometry onto $B_d (s)$.
\end{proof}

\begin{CorRef}{ChIso}
Given a finitely generated $R[\Gamma]$-module $F$ and two choices
of finite generating sets $\Sigma_1$ and $\Sigma_2$, the filtered
$R$-modules $s (F, \Sigma_1)$ and $s (F, \Sigma_2)$ are isomorphic
as $\Gamma$-filtered $R$-modules.
\end{CorRef}

\begin{proof}
The identity map and its inverse are boundedly controlled as maps
between $s (F, \Sigma_1)$ and $s (F, \Sigma_2)$ by Lemma
\refL{AlmPull}.
\end{proof}

\begin{CorRef}{FilProps}
Finitely generated $R[\Gamma]$-modules $F$ with filtrations $s (F,
\Sigma)$, with respect to arbitrary finite generating sets
$\Sigma$, are locally finitely generated and lean.  If $s (F,
\Sigma)$ is insular and $\Sigma'$ is another finite generating set
then $s (F, \Sigma')$ is also insular.
\end{CorRef}

\begin{proof}
For a finite subset $S$, the submodule $F(S)$ is generated by the
finite set $S\Sigma$. Since $F(x) = \langle x \Sigma \rangle_R$,
\[
F(S) = \sum_{x \in \Sigma} \langle x \Sigma \rangle_R = \langle S
\Sigma \rangle_R,
\]
so $s(F,\Sigma)$ is $0$-lean.  The second
claim follows from Corollary \refC{ChIso}.
\end{proof}

\begin{DefRef}{BR[G]}
Let $\B (R[\Gamma])$ be the full subcategory of $\Modf
(R[\Gamma])$ on $R$-modules $F$ which are
\textit{strict} as filtered modules $s (F, \Sigma)$ with respect
to some choice of the finite generating set $\Sigma$.

Let $\Btimes (R[\Gamma])$ be the category of objects which are
pairs $(F, \Sigma)$ with $F$ in $\B (R[\Gamma])$ and
$\Sigma$ a finite generating set for $F$. The morphisms are
the $R[\Gamma]$-homomorphisms between the modules.
\end{DefRef}

Lemma \refL{AlmPull} shows that the map
\[
s \colon \Btimes (R[\Gamma]) \longrightarrow \B (\Gamma,R)
\]
described in Definition \refD{AndBack} is a functor. In fact, it
is a functor
\[
s_{\Gamma} \colon \Btimes (R[\Gamma]) \longrightarrow
\B_{\Gamma,0} (\Gamma,R)^{\Gamma}
\]
by interpreting $s_{\Gamma} (F, \Sigma) = (F, \psi)$ with $F = s (F,
\Sigma)$ and $\psi (\gamma) \colon F \to \gamma F$ induced from $s
\sigma \mapsto \gamma^{-1} s \sigma$. Since
\[
(\gamma F)(S) =
\langle \gamma^{-1} (S) \Sigma \rangle_R,
\]
it follows that the
object $s_{\Gamma} (F, \Sigma)$ lands in $\B_{\Gamma,0}
(\Gamma,R)^{\Gamma}$, and $s$ sends all $R[\Gamma]$-homomorphisms
to $\Gamma$-equivariant homomorphisms.

\begin{LemRef}{ForAppr}
Let $F \in \B_{\Gamma,0} (\Gamma, R)^{\Gamma}$ and let $\Sigma$ be
a finite generating set for the $R[\Gamma]$-module $F$. Then the
identity homomorphism
\[
\id \colon s_{\Gamma} (F, \Sigma) \longrightarrow F
\]
is
boundedly controlled with respect to the induced and the original
filtrations of $F$.
\end{LemRef}

\begin{proof}
If $\Sigma$ is contained in $F(B_d (e))$, where $e$ is the
identity element in $\Gamma$, then
\[
\gamma \Sigma \subset F(B_d
(\gamma))
\]
for all $\gamma \in \Gamma$, and
\[
s (F, \Sigma) (S) =
\langle S \Sigma \rangle_R \subset F(S[d])
\]
for all subsets $S
\subset \Gamma$.
\end{proof}

Both functors $s$ and $s_{\Gamma}$ are additive with respect to
the additive structure in $\Btimes (R[\Gamma])$ where the biproduct is given by
\[
(F, \Sigma_F) \oplus (G, \Sigma_G) = (F \oplus G, \Sigma_F \times
\Sigma_G).
\]
Let the \textit{admissible monomorphisms} $\phi \colon
(F, \Sigma_F) \to (G, \Sigma_G)$ in $\Btimes (R[\Gamma])$ be the
injections $\phi \colon F \to G$ of $R[\Gamma]$-modules $\phi$
such that
\[
s (\phi) \colon s (F, \Sigma_F) \longrightarrow s (G, \Sigma_G)
\]
is
a boundedly bicontrolled homomorphism of $\Gamma$-filtered
$R$-modules. This is equivalent to requiring that $s(\phi)$ be an
admissible monomorphism in $\B (\Gamma,R)$. Let the
\textit{admissible epimorphisms} be the morphisms $\phi$ such that
$s(\phi)$ are admissible epimorphisms in $\B (\Gamma,R)$.

\begin{PropRef}{ItIsEx}
The choice of admissible morphisms defines an exact structure on
$\Btimes (R[\Gamma])$ such that both $s$ and $s_{\Gamma}$ are
exact functors.
\end{PropRef}

\begin{proof}
When checking Quillen's axioms in $\Btimes (R[\Gamma])$, all
required universal constructions are performed in $\B (R[\Gamma])$
with the canonical choices of finite generating sets. In
particular, $\Sigma$ in the pushout $B \cup_A C$ is the image of
the product set $\Sigma_B \times \Sigma_C$ in $B \times C$.

The
fact that all admissible
morphisms are boundedly bicontrolled in $\B (\Gamma, R)$ or
$\B_{\Gamma,0} (\Gamma, R)^{\Gamma}$ follows from the proof of
Theorem \refT{UbisWAb}. Exactness of $s$ and $s_{\Gamma}$ is
immediate.
\end{proof}

\begin{DefRef}{ghjjk}
We give $\B (R[\Gamma])$ the minimal exact structure that makes
the forgetful functor
\[
p \colon \Btimes (R[\Gamma]) \longrightarrow \B (R[\Gamma])
\]
sending $(F, \Sigma)$ to $F$ an exact functor.

In other words, an $R[\Gamma]$-homomorphism
$\phi \colon F \to G$ is an \textit{admissible monomorphism} or
\textit{epimorphism} if for some choice of finite generating sets,
\[
\phi \colon (F, \Sigma_F) \longrightarrow (G, \Sigma_G)
\]
is respectively an
admissible monomorphism or epimorphism in $\Btimes (R[\Gamma])$.

Corollary \refC{ChIso} shows that if $\phi \colon F \to G$ is
boundedly bicontrolled as a map of filtered $R$-modules $s (F,
\Sigma_F) \to s (G, \Sigma_G)$ then it is boundedly bicontrolled
with respect to any other choice of finite generating sets, so
this structure is well-defined.
\end{DefRef}

\begin{NotRef}{NewExStr}
The new exact category will be referred to as $\B (R[\Gamma])$,
with the corresponding $K$-theory spectrum $\Gnc (R[\Gamma])$.
\end{NotRef}

Let $(F,\psi)$ be an object of $\B_{\Gamma,0}
(\Gamma,R)^{\Gamma}$. One may think of $\gamma F \in \B
(\Gamma,R)$, $\gamma \in \Gamma$, as the module $F$ with a new
$\Gamma$-filtration. Now the $R$-module structure $\eta \colon R
\to \End F$ induces an $R[\Gamma]$-module structure
\[
\eta (\psi)
\colon R[\Gamma] \longrightarrow \End F
\]
given by
\[
\sum_{\gamma} r_{\gamma} \gamma \mapsto \sum_{\gamma} \eta
(r_{\gamma})
 \psi (\gamma)
\]
since the sums are taken over a finite subset of $\Gamma$. It is
easy to see that this defines a map
\[
\pi \colon \B_{\Gamma,0} (\Gamma,R)^{\Gamma} \longrightarrow \B
(R[\Gamma])
\]
by sending $(F, \psi)$ to $F$, so that $p = \pi s_{\Gamma}$.
Notice however that in general $\pi$ is not exact as the identity
homomorphism in Lemma \refL{ForAppr} is not necessarily an
isomorphism.

In the rest of the paper we assume $\Gamma$ is torsion-free. The
exact functors $p$ and $s_{\Gamma}$ induce maps in nonconnective
$K$-theory
\[
\Gnc (R[\Gamma]) \xleftarrow{\ p} \Knc \Btimes (R[\Gamma])
\xrightarrow{\ s_{\Gamma}\ } \GncGammazero (\Gamma,R)^{\Gamma}.
\]
We claim that both of these maps are weak equivalences.

\begin{PropRef}{Debase}
The functor $f$ induces a weak equivalence
\[
\Knc \Btimes (R[\Gamma]) \simeq \Gnc (R[\Gamma]).
\]
\end{PropRef}

\begin{proof}
This follows from the Approximation theorem applied to $p'$. The
two categories are saturated, and $\Btimes (R[\Gamma])'$ has a
cylinder functor satisfying the cylinder axiom which is
constructed as the canonical homotopy pushout with the canonical
product basis, see section 1 of \cite{rTtT:90}.

The first
condition of the Approximation theorem is clear. For the second
condition, let $(F^{\subdot}_1, \Sigma^{\subdot}_1)$ be a complex
in $\Btimes (R[\Gamma])$ and let $g \colon F^{\subdot}_1 \to
F^{\subdot}_2$ be a chain map in $\B (R[\Gamma])'$. For each
$R[\Gamma]$-module $F^i_2$ choose any finite generating set
$\Sigma^i_2$, then using $f=g$ and $g'=\id$, we have $g = g'
p(f)$.
\end{proof}

\begin{PropRef}{InjGR}
The functor $s_{\Gamma}$ induces a weak homotopy equivalence
\[
\Knc \Btimes (R[\Gamma]) \simeq \GncGammazero (\Gamma,
R)^{\Gamma}.
\]
\end{PropRef}

\begin{proof}
The target category is again saturated and has a cylinder functor
satisfying the cylinder axiom. To check condition 2 of the
approximation theorem, let
\[
E^{\subdot} \colon \quad 0 \longrightarrow (E^1, \Sigma_1)
\longrightarrow (E^2, \Sigma_2) \longrightarrow\ \dots\
\longrightarrow (E^n, \Sigma_n) \longrightarrow 0
\]
be a complex in $\Btimes (R[\Gamma])$,
\[
(F^{\subdot},\psi.) \colon \quad 0 \longrightarrow (F^1, \psi_1)
\xrightarrow{\ f_1\ } (F^2, \psi_2) \xrightarrow{\ f_2\ }\ \dots\
\xrightarrow{\ f_{n-1}\ } (F^n, \psi_n) \longrightarrow 0
\]
be a complex in $\B_{\Gamma,0} (\Gamma,R)^{\Gamma}$, and
\[
g \colon
s'_{\Gamma} (E^{\subdot}) \longrightarrow (F^{\subdot}, \psi.)
\]
be a chain
map. Each $F^i$ can be thought of as an $R[\Gamma]$-module, and
there is a chain complex
\[
F^{\subdot} \colon \quad 0 \longrightarrow F^1 \xrightarrow{\ f_1\
} F^2 \xrightarrow{\ f_2\ }\ \dots\ \xrightarrow{\ f_{n-1}\ } F^n
\longrightarrow 0
\]
in $\Modf (R[\Gamma])$. Choose arbitrary finite generating sets
$\Omega_i$ in $F^i$ for all $1 \le i \le n$. Now
\[
\pi_{\Omega} F^{\subdot} \colon \quad 0 \longrightarrow (F^1,
\Omega_1) \xrightarrow{\ f_1\ } (F^2, \Omega_2) \xrightarrow{\
f_2\ }\ \dots\ \xrightarrow{\ f_{n-1}\ } (F^n, \Omega_n)
\longrightarrow 0
\]
is a chain complex in $\Btimes (R[\Gamma])$.

The chain map $g$
is degree-wise an $R[\Gamma]$-homomorphism, so there is a
corresponding chain map
\[
f \colon E^{\subdot} \longrightarrow \pi_{\Omega}
F^{\subdot}
\]
which coincides with $g$ on modules. On the other
hand, the degree-wise identity gives a chain map
\[
g' \colon
s'_{\Gamma} (\pi_{\Omega} F^{\subdot}) \longrightarrow F^{\subdot}
\]
in
$\B_{\Gamma,0} (\Gamma, R)^{\Gamma}$ by Lemma \refL{ForAppr}. This
$g'$ is a quasi-isomorphism, as required.
\end{proof}

\begin{CorRef}{InjGRreal}
Let $\Gamma$ be a finitely generated torsion-free group and $R$ be
a Noetherian ring. There is a weak equivalence
\[
\GncGammazero (\Gamma ,R)^{\Gamma} \simeq \Gnc (R[\Gamma]).
\]
\end{CorRef}

\begin{CorRef}{InjGR}
Let $\Gamma$ be a finitely generated torsion-free group acting
freely, properly discontinuously and cocompactly on a proper
metric space $X$ and let $R$ be a Noetherian ring.  There is a
weak homotopy equivalence
\[
\GncGammazero (X ,R)^{\Gamma} \simeq \Gnc (R[\Gamma]).
\]
\end{CorRef}

\begin{proof}
Let $p \colon X \to \mathrm{point}$ be the geometric collapse. For
any $x \in X$ such that the embedding $i$ of the orbit $\Gamma x$
with the word metric is commensurable in $X$, there is a
commutative diagram
\[
\begin{CD}
\GncGammazero (\Gamma x,R)^{\Gamma}  @>{\pi_{\ast}}>> \Gnc (R[\Gamma]) \\
@VV{i_{\ast}}V @V{=}VV\\
\GncGammazero (X,R)^{\Gamma}  @>{\pi_{\ast}}>> \Gnc (R[\Gamma])
\end{CD}
\]
The top $\pi_{\ast}$ is a weak equivalence by Corollary
\refC{InjGRreal}. The vertical map $i_{\ast}$ is a weak equivalence as
in Proposition \refP{IsomOK}, so the lower map $\pi_{\ast}$
is a weak equivalence.
\end{proof}

For a discrete group $\Gamma$ and a ring $R$ there is an assembly
map
\[
A_K \colon B\Gamma_{+} \wedge \Knc (R) \longrightarrow \Knc
(R[\Gamma]).
\]
When the ring $R$ is regular
Noetherian, for example the integers $\mathbb{Z}$, the spectra $\Gnc (R)$ and $\Knc (R)$ and, therefore, $B\Gamma_{+} \wedge \Gnc (R)$ and $B\Gamma_{+} \wedge \Knc (R)$ can be naturally
identified.

\begin{DefRef}{AssG}
Let $\Gamma$ be a finitely generated group and $R$ be  a regular Noetherian ring.
The \textit{assembly map in $G$-theory}
\[
A_G \colon B\Gamma_{+} \wedge \Gnc (R) \longrightarrow \Gnc (R[\Gamma])
\]
is the
composition of $A_K$ and the canonical Cartan map
\[
C \colon \Knc
(R[\Gamma]) \longrightarrow \Gnc (R[\Gamma])
\]
induced by inclusion of
categories.

The \define{integral Novikov conjecture} in algebraic
$G$-theory is the statement that this is a split injection of
spectra.
\end{DefRef}

\begin{RemRef}{GtoK}
Notice that whenever the assembly map $A_G$ is split
injective, the map $A_K$ is also split injective, so
this conjecture is stronger than the $K$-theoretic conjecture when
the ring $R$ is regular.
\end{RemRef}

\begin{RemRef}{Luck}
The standard exact structure on $\Modf (R[\Gamma])$ has all injective and surjective $R[\Gamma]$-homomorphisms with
finitely generated cokernels and kernels as admissible morphisms
so that the exact sequences are the traditional short exact
sequences. Let the corresponding
$K$-theory spectrum be $\Gncm (R[\Gamma])$.
One might attempt to
replace $\Gnctimes (R[\Gamma])$ with $\Gncm (R[\Gamma])$ as the target of the assembly $A_G$.
However, W.~L\"{u}ck has pointed out that this map would not be
weakly injective even in the simple case when $R$ is a commutative
ring and $\Gamma$ is the free group on two generators, cf.~Remark
2.23 in \cite{wL:98}.

This underscores the importance of choosing the coarse version $\Gnc (R[\Gamma])$ as our approximation of $\Knc (R[\Gamma])$.
\end{RemRef}

The equivariant \textit{assembly map in boundedly controlled $G$-theory} can be defined as in
\cite{gC:95}.
For any Noetherian ring $R$, this is an equivariant natural transformation
\[
\alpha_G \colon \hlf (X, \Gnc (R)) \longrightarrow \GncGammazero (X,R)
\]
from the equivariant locally
finite homology $\hlf (X, \Gnc (R))$ to the equivariant bounded $G$-theory $\GncGammazero (X,R)$, see
Definition II.14, loc.~cit.
If $\Gamma$ is torsion-free and acts freely cocompactly on $X$, one also has weak equivalences
\[
\hlf (X, \Gnc (R))^{\Gamma} \simeq B\Gamma_{+} \wedge \Gnc (R).
\]
The fixed point spectra and the induced maps fit in the commutative diagram
\begin{equation}
\begin{CD}
B\Gamma_{+} \wedge \Gnc (R) @>{\alpha_K^{\Gamma}}>> \KncGammazero
(X,R)^{\Gamma} @>{\simeq}>>
\Knc (R[\Gamma])\\
@V{\alpha_G^{\Gamma}}VV @VVV @VV{C}V \\
\GncGammazero (X,R)^{\Gamma} @<{\simeq}<{s_{\Gamma \ast}}< \Knc
\Btimes (R[\Gamma]) @>{=}>> \Gnc
(R[\Gamma])
\end{CD}
\tag{$\dagger \dagger$}
\end{equation}

If we restrict our attention to the
equivariant objects in $\KncGammazero (X,R)^{\Gamma}$ that have
the parametrization function map the generating set $B$ to a single point
then the map
\[
i \colon \KncGammazero (X,R)^{\Gamma} \longrightarrow
\GncGammazero (X,R)^{\Gamma}
\]
is well-defined and fits as the diagonal in
the left square.

\begin{RemRef}{HUMKIU}
When $R$ is a regular Noetherian ring,
the assembly map in boundedly controlled $G$-theory
\[
\alpha_G \colon \hlf (X, \Gnc (R)) \longrightarrow \GncGammazero (X,R)
\]
can be identified up to homotopy with the composition
\[
\hlf (X, \Knc (R))\xrightarrow{\ \alpha_K\ }
\KncGammazero (X,R) \xrightarrow{\ i\ } \GncGammazero (X,R).
\]
Now there is a homotopy commutative square
\[
\begin{CD}
B\Gamma_{+} \wedge \Knc (R) @>{A_G}>> \Gnc
(R[\Gamma]) \\
@V{\simeq}VV @VV{\simeq}V \\
\hlf (X, \Gnc (R))^{\Gamma} @>{\alpha_G^{\Gamma}}>> \GncGammazero (X,R)^{\Gamma}
\end{CD}
\]
\end{RemRef}

Controlled excision is the main technical tool from controlled
$K$-theory used to prove integral Novikov conjectures. It is used
to see that in specific cases the equivariant assembly map in bounded
$K$-theory is a homotopy equivalence. In particular, the argument
in \cite{gCbG:04} applies to groups of finite asymptotic dimension
which have a finite classifying space.

\begin{ThmRef}{NCGTH}
The fixed point map of spectra
\[
\alpha_G^{\Gamma} \colon \hlf (X, \Gnc (R))^{\Gamma} \longrightarrow \GncGammazero (X,R)^{\Gamma}
\]
is a split injection for any geometrically finite group $\Gamma$
of finite asymptotic dimension and a Noetherian ring $R$.
\end{ThmRef}

\begin{proof}
The main step in the proofs of the Novikov conjecture in algebraic
$K$-theory \cite{gC:95,gCbG:04} to which we referred above is the
application of homotopy fixed points to reduce the study of the
map $A_K$ to the nonequivariant study of the equivariant map
$\alpha_K$. This is shown to be a weak equivalence by using
controlled excision to compute the target. With the excision
results from section 3 and the equivariant properties established
here, the proofs can be repeated
verbatim obtaining splittings of the assembly maps $\alpha_G$ for
the same collection of groups.
\end{proof}

\begin{CorRefName}{CorNCGTH}{Novikov conjecture in $G$-theory}
The $G$-theoretic assembly map
\[
A_G \colon B\Gamma_{+} \wedge \Gnc (R) \longrightarrow \Gnc
(R[\Gamma])
\]
is a split injection for any geometrically finite group $\Gamma$
of finite asymptotic dimension and a regular Noetherian ring $R$.
\end{CorRefName}

\begin{proof}
If $R$ is regular Noetherian, the two maps $A_G$ and $\alpha_G^{\Gamma}$
can be identified as in Remark \refR{HUMKIU}.
\end{proof}


\begin{thebibliography}{99}

\bibitem{dAhM:88}
{D.R.~Anderson and H.J.~Munkholm}, \textit{Boundedly controlled
topology}, Lecture Notes in Mathematics {\bf 1323},
Springer-Verlag (1988).

\bibitem{hB:68}
{H.~Bass}, \textit{Algebraic $K$-theory}, W.~A.~Benjamin, Inc.
(1968).

\bibitem{mCeP:97}
{M.~Cardenas and E.K.~Pedersen}, {\it On the Karoubi filtration of
a category}, $K$-theory, {\bf 12} (1997), 165--191.

\bibitem{gC:95}
{G. Carlsson}, {\it Bounded $K$-theory and the assembly map in
algebraic $K$-theory}, in {\it Novikov conjectures, index theory
and rigidity}, {\it Vol. 2} (S.C.~Ferry, A.~Ranicki, and
J.~Rosenberg, eds.), Cambridge U. Press (1995), 5--127.

\bibitem{gCbG:03}
{G.~Carlsson and B.~Goldfarb}, \textit{On homological coherence of
discrete groups}, J.~Algebra {\bf 276} (2004), 502--514.

\bibitem{gCbG:04}
\bysame, \textit{The integral K-theoretic Novikov conjecture for
groups with finite asymptotic dimension}, Inventionnes Math. {\bf
157} (2004), 405--418.

\bibitem{gCbG:06}
\bysame, \textit{Algebraic $K$-theory of geometric groups}, in
preparation.

\bibitem{dF:80}
{D.R.~Farkas}, \textit{Group rings: an annotated questionnaire},
Comm. Algebra \textbf{8} (1980), 585--602.

\bibitem{pGmZ:67}
{P. Gabriel and M. Zisman}, \textit{Calculus of fractions and
homotopy theory}, Springer-Verlag (1967).

\bibitem{aiG:96}
{A.I.~Generalov}, \textit{Relative homological algebra. Cohomology of Categories, posets and coalgebras}, in
Handbook of Algebra, Vol.~1 (M.~Hazewinkel, ed.), 1996, Elsevier
Science, 611--638.

\bibitem{iHeP:99}
{I.~Hambleton and E.K.~Pedersen}, \textit{Compactifying infinite
group actions}, Contemp. Math. {\bf 258} (2000), 203--212.

\bibitem{rH:77}
{R. Hartshorne}, {\it Algebraic geometry}, Springer-Verlag (1977).

\bibitem{bK:90}
{B.~Keller}, \textit{Chain complexes and stable categories},
Manuscripta Math. \textbf{67} (1990), 379--417.

\bibitem{bK:96}
\bysame, \textit{Derived categories and their uses}, in
Handbook of Algebra, Vol.~1 (M.~Hazewinkel, ed.), 1996, Elsevier
Science, 671--701.

\bibitem{wL:98}
{W.~L\"{u}ck}, \textit{Dimension theory of arbitrary modules over
finite von Neumann algebras and $L^2$-Betti numbers II:
Applications to Grothendieck groups, $L^2$-Euler characteristics
and Burnside groups}, J. reine angew. Math. \textbf{496} (1998),
213--236.

\bibitem{sM:71}
{S.~Mac Lane}, \textit{Categories for the working mathematician},
Springer-Verlag (1971).

\bibitem{ePcW:85}
{E.K. Pedersen and C. Weibel}, {\it A nonconnective delooping of
algebraic $K$-theory}, in {\it Algebraic and geometric topology}
(A.~Ranicki, N.~Levitt, and F.~Quinn, eds.), Lecture Notes in
Mathematics {\bf 1126}, Springer-Verlag (1985), 166--181.

\bibitem{ePcW:89}
{\bysame}, {\it $K$-theory homology of spaces}, in {\it Algebraic
topology} (G.~Carlsson, R.L.~Cohen, H.R.~Miller, and
D.C.~Ravenel, eds.), Lecture Notes in Mathematics {\bf 1370},
Springer-Verlag (1989), 346--361.

\bibitem{ePcW:97}
{\bysame}, unpublished.

\bibitem{dQ:73}
{D. Quillen}, {\it Higher algebraic $K$-theory: I}, in {\it
Algebraic $K$-theory I} (H.~Bass, ed.), Lecture Notes in
Mathematics {\bf 341}, Springer-Verlag (1973), 77--139.

\bibitem{mS:99}
{M.~Schlichting}, \textit{Delooping the $K$-theory of exact
categories and negative $K$-groups}, U. Paris VII thesis (2000).

\bibitem{mS:03}
{\bysame}, \textit{Delooping the $K$-theory of exact categories},
Topology {\bf 43} (2004), 1089--1103.

\bibitem{rS:89}
{R.E.~Staffeldt}, \textit{On the fundamental theorems of algebraic
$K$-theory}, $K$-theory {\bf 1} (1989), 511--532.

\bibitem{rT:83}
{R.W. Thomason}, {\it The homotopy limit problem}, in {\it
Proceedings of the Northwestern homotopy theory conference} (H.R.
Miller and S.B. Priddy, eds.), Cont. Math. {\bf 19} (1983),
407--420.

\bibitem{rTtT:90}
{R.W. Thomason and Thomas Trobaugh}, {\it Higher algebraic
$K$-theory of schemes and of derived categories}, in {\it The
Grothendieck Festschrift}, {\it Vol. III}, Progress in Mathematics
{\bf 88}, Birkh\"{a}user (1990), 247--435.

\bibitem{fW:85}
{F. Waldhausen}, {\it Algebraic $K$-theory of spaces}, in {\it
Algebraic and geometric topology} (A. Ranicki, N. Levitt, and F.
Quinn, eds.), Lecture Notes in Mathematics {\bf 1126},
Springer-Verlag (1985), 318--419.

\end{thebibliography}
\end{document}